\documentclass[11pt]{article}
\usepackage{amssymb}
\usepackage{mathrsfs}
\usepackage{latexsym,amsmath,amssymb,amsfonts,epsfig,graphicx,cite,psfrag}
\usepackage{eepic,color,colordvi,amscd}
\usepackage{ebezier}
\usepackage{verbatim}
\usepackage{subfigure}
\usepackage{enumitem}
\usepackage{algorithm}
\usepackage{algorithmic}
\usepackage{hyperref}
\usepackage{amsthm}
\usepackage{longtable}
\usepackage{comment}
\usepackage{color}
\usepackage{graphicx} 
\usepackage{epstopdf}
\definecolor{blue}{rgb}{0,0,0.9}
\definecolor{red}{rgb}{0.9,0,0}
\definecolor{green}{rgb}{0,0.9,0}

\theoremstyle{plain}
\newtheorem{remark}{Remark}

\newtheorem{assumption}{Assumption}
\newtheorem{definition}{Definition}
\newtheorem{proposition}{Proposition}
\newtheorem{theorem}{Theorem}
\newtheorem{claim}{Claim}
\newtheorem{lemma}{Lemma}
\usepackage{longtable}

\def\Tr{\mathrm{Tr}}
\def\H{\mathcal{H}}
\def\<{\big\langle}
\def\>{\big\rangle}
\def\E{\mathcal{E}}
\def\M{\mathcal{M}}
\def\T{{\rm T}}
\def\D{{\rm D}}
\def\DD{{\rm Diag}}
\def\A{\mathcal{A}}

\def\K{\mathcal{K}}
\def\L{\mathcal{L}}
\def\F{\mathcal{F}}
\def\I{\mathcal{I}}
\def\J{\mathcal{J}}
\def\R{\mathbb{R}}

\def\Re{{\rm Rtr}}

\def\g{{\rm grad}}
\def\G{\mathcal{G}}
\def\N{\mathcal{N}}

\def\S{\mathbb{S}}

\def\rr{{\rm rank}}
\def\dist{{\rm dist}}
\def\dd{{\rm diag}}
\def\P{{\rm Proj}}

\def\BA{\Omega_{r}}
\def\O{\mathcal{O}}

\def\({\left(}
\def\){\right)}
\def\hR{\widehat{R}}

\def\hS{\widehat{S}}

\def\tA{\tilde{\A}}

\def\fval{{\rm fval}}
\let\svthefootnote\thefootnote
\newcommand\blankfootnote[1]{%
	\let\thefootnote\relax\footnotetext{#1}%
	\let\thefootnote\svthefootnote%
}

\usepackage{amsopn}

\begin{document}

	\title{A feasible method for general convex low-rank SDP problems}
		
	\author{Tianyun Tang
\thanks{Department of Mathematics, National
         University of Singapore, Singapore
         119076 ({\tt ttang@u.nus.edu}).
         }, \quad 
	 Kim-Chuan Toh\thanks{Department of Mathematics, and Institute of 
Operations Research and Analytics, National
         University of Singapore, 
       Singapore
         119076 ({\tt mattohkc@nus.edu.sg}).  The research of this author is supported
by the Ministry of Education, Singapore, under its Academic Research Fund Tier 3 grant call (MOE-2019-T3-1-010).}
	 }
	
	\date{\today}
	
\maketitle

\begin{abstract}
In this work, we consider the low rank decomposition (SDPR) of general convex semidefinite programming problems (SDP) that contain both a positive semidefinite matrix and a nonnegative vector as variables. We develop a rank-support-adaptive feasible method to solve (SDPR) based on Riemannian optimization. The method is able to escape from a saddle point to ensure its convergence to a 
global optimal solution for generic constraint vectors.
We prove its global convergence and local linear convergence without assuming that the objective function is twice differentiable. Due to the special structure of the low-rank SDP problem, our algorithm can achieve better iteration complexity than existing results 
for more general smooth nonconvex problems. In order to overcome the degeneracy issues of SDP problems, we develop two strategies based on random perturbation and dual refinement. These techniques enable us to solve some primal degenerate SDP problems efficiently, for example, Lov\'{a}sz theta SDPs. Our work is a step forward in extending the application range of Riemannian optimization approaches for solving SDP problems. Numerical experiments are conducted to verify the efficiency and robustness of our method.
\end{abstract}

\bigskip
\noindent{\bf keywords:} semidefinte programming, feasible method, Riemannian optimization
\\[5pt]

\section{Introduction}
\subsection{Low rank SDP problem}\label{LRSDP}
The standard  linear SDP problem is as follows:
\vskip-4mm
\begin{equation}\label{SDP0}
\min\left\{ \<C,X\>:\ \A(X)=b,\ X\in \S^n_+ \right\},
\end{equation}
where $C\in \S^n,$ $b\in \R^m$ and $\A:\S^n\rightarrow \R^m$ is a linear mapping. Problem \eqref{SDP0} can be solved with guaranteed convergence by many well-developed convex solvers like SDPT3 \cite{TTT,T3Q}, MOSEK \cite{mosek} and SDPNAL \cite{SDPNALp2,SDPNALp1,SDPNAL}. However, due to the dimensionality of $O(n^2)$
of the matrix variable,  convex solvers are inefficient when $n$ is large. To overcome this issue, Burer and Monteiro \cite{BM3,BM1} proposed the following low rank decomposition of \eqref{SDP0}, which we also call as the factorized model:
\begin{equation}\label{SDPLR}
 \min\left\{ \<C,RR^\top\>:\ \A(RR^\top)=b,\ R\in \R^{n\times r} \right\}.
\end{equation}
From the well known low rank property of problem \eqref{SDP0} \cite{BAI,FSL,PG}, problem \eqref{SDP0} and \eqref{SDPLR} are equivalent when $r(r+1)/2\geq m.$ By choosing $r=\lceil \sqrt{2m}\rceil,$ the dimension of \eqref{SDPLR} becomes $O(nm^{1/2}),$ which is much smaller than $O(n^2)$ if $m=o(n^2).$  


\subsection{Feasible method and its limitation}\label{Limit}
There are mainly two types of methods to solve \eqref{SDPLR}. One is augmented Lagrangian method (ALM). It has been implemented as a general linear SDP solver by Burer and Monteiro. The other one is feasible method\footnote{Throughout this paper, by a feasible method, we always mean Riemannian optimization approach.} based on Riemannian optimization (see \cite{syn,Staircase1,Staircase2,WYmani}). This method assumes that the feasible set of \eqref{SDPLR} is a Riemannian manifold.
Thus, one can use Riemannian optimization method such as Riemannian gradient descent and Riemannian trust region method to solve the problem \cite{manibook,Intromani,manopt,Geoopt,CDopt1,CDopt}. Apart from ALM and Riemannian optimization methods, Bellavia et al. \cite{bellavia2021relaxed} developed a rank-adaptive interior point method for low-rank linear SDP problems and applied it to matrix completion problems.

It is well-known that the efficiency of ALM relies heavily on the penalty parameter and the termination criterion for solving its subproblems. 
Compared with ALM, feasible method is penalty-free and one can directly decrease the objective function value on the manifold just like unconstrained optimization. 
Another advantage of feasible method is that it can find a solution of acceptable accuracy after only a few iterations while the primal feasibility is maintained. 
However, the application range of feasible method is quite narrow. As far as we know, feasible method has never been used in a general SDP solver like interior point method and ALM. 
Current feasible method is mainly designed for SDP problems with block-diagonal equality constraints such as the SDP relaxations of max-cut problems, robust PCA and synchronization of rotations \cite{Staircase1,erdogdu2022convergence}. 
The feasible sets of such SDP problems are simple smooth manifolds like oblique manifold and Stiefel manifold. 
Different from these special problems, the constraints of a general SDP problem may be more complicated with inequality constraints and its feasible set may contain singular\footnote{Throughout this paper, ``singular" means that the linearly independent constraint qualification (LICQ) is not satisfied.} points. 
The worst case is that the optimal solution lies exactly on a singular point. This will cause numerical issues as the iteration points approach it. 
Recently, Wang et al. \cite{wang2023solving} and Wen et al. \cite{wang2021decomposition} developed Riemannian ALMs that solve \eqref{SDPLR} by separating its constraints into two parts such that one of them is a Riemannian manifold. 
These methods are not feasible methods and the manifolds they considered are still the special manifolds just mentioned above. In \cite{GEP,SQK}, we have applied a feasible method with rounding and escaping strategy to solve some special SDP problems with finite number of singular points. 
However, for a general SDP problem, there might be infinite number (even uncountable) of singular points. It would be challenging to extend that idea to general SDP problems. 

In this paper, our goal is to demonstrate that feasible methods can also be developed to solve a general SDP problem. However, we should emphasize that unlike other previous works
on Riemannian optimization approaches for solving an SDP with special constraints, we target SDP problems with 
general linear constraints and the corresponding feasible sets of the factorized problems 
\eqref{SDPLR} are algebraic varieties rather
than special manifolds. Thus, even though our rank-support-adaptive feasible method is motivated
by Riemannian optimization approaches, we must overcome the inherent theoretical and numerical difficulties brought about by singular points in the feasible sets. In these aspects,
our design of the feasible method for solving \eqref{SDPLR}, and its theoretical analysis and
numerical implementation, are much more challenging than many existing works which target
only SDP problems with special constraints. We will state our contributions and their relationship to existing works in the literature in the following subsections. 

\subsection{A general SDP model and its low rank decomposition}\label{Contri}
We consider the following convex (possibly nonlinear) SDP problem:
\vskip-4mm
\begin{equation}\label{SDP}
({\rm SDP})\ \min\left\{ \phi(X,x):\ \K(X,x)=b,\ X=\begin{pmatrix}I_k& * \\ * & *\end{pmatrix}\in \S_+^n,\ x\in \R^p_+\right\}.
\end{equation}
In (SDP), $\phi:\S^n\times \R^p\rightarrow \R$ is a convex continuously differentiable function, $b\in \R^m$, $\K:\S^n\times \R^p\rightarrow \R^m$ is a  linear mapping such that for any $(X,x)\in \S^n\times \R^p,$ 
\begin{equation}\label{Kmap}
\K\(X,x\):=\A(X)+Bx,
\end{equation}
where $B\in \R^{m\times p}$ and $\A:\S^n\rightarrow \R^m$ is a linear mapping such that 
$$\A(X):=\(\<A_1,X\>,\< A_2,X\>,\ldots,\<A_m,X\>\)^\top$$
 for data matrices $A_1,A_2,\ldots,A_m\in \S^n.$ (SDP) is more general than problem \eqref{SDP0} because the objective function may not be linear and there is an additional nonnegative vector variable $x\in \R^{p}_+.$ This vector variable may occur in SDP problems with inequality constraints, where it is used as the slack variable. Note that 
  (SDP) has an implicit fixed constraint given by $X_{1:k,1:k}=I_k,$ but the case where $k=0$ is allowed.
 This constraint, with $k=1$, appears in the SDP relaxations of quadratically constrained quadratic programming with linear terms \cite{QCQP1}. We 
 intentionally separate the fixed constraints from $\A\(X\)=b$ because of its special structure. Now, we consider the following low rank decomposition of (SDP):
\begin{equation}\label{SDPR}
({\rm SDPR})\ \min\left\{ \phi\big(\hR\hR^\top,y\circ y\big):\ \K\big(\hR\hR^\top,y\circ y\big) = b, \hR=[I_{k,r};R]\in \R^{n\times r},\ y\in \R^p\notag\right\}, 
\end{equation}
where $k\leq r\leq n$, ``$\circ$'' denotes the Hadamard product\footnote{If $x$ is the slack variable of some inequality constraints, then the vector $y\circ y$ is called squared slack variable.} and
$I_{k,r} = [I_k, 0_{k\times (r-k)}]$. In (SDPR), we consider $R\in \R^{(n-k)\times r}$ and $y\in \R^p$ as unknown variables. Different from the traditional Burer-Monteiro factorization, we fixed the first $k$ rows of the factorized matrix to be the $k$ by $r$ identity matrix $I_{k,r}.$ By doing so, the constraints $X_{1:k,1:k}=I_k$ is automatically satisfied. This idea has also been used in \cite{SQK} to deal with the constraint $X_{11}=1.$ In (SDPR), our constraint is more general. It is easy to see that the optimal value of (SDPR) is an upper bound of (SDP). Moreover, when (SDP) has an optimal solution $\(X,x\)$ such that $\rr(X)= r$, $X$ can be factorized as $\hR\hR^\top$ with $\hR\in \R^{n\times r}$ and $\hR_{1:k,:}=I_{k,r}.$ In this case, $\big(R:=\hR_{k+1:n,:},y:=\sqrt{x}\big)$\footnote{Here we use elementwise square root.} is feasible for (SDPR) and so (SDPR) is equivalent to (SDP). 

\subsection{Optimality condition for (SDPR)}\label{Subsec-OPTcon}
There are many works in the literature studying the optimality condition of low rank SDP problems (see \cite{Boumal1,Boumal2,Boumal3,BM2,CIFrank,Boumal4,GEP}). They prove that for a linear SDP problem with a generic cost matrix $C,$ a second order stationary point of the factorized SDP problem is also a global optimal solution provided that certain constraint nondegeneracy conditions and rank conditions hold\footnote{The rank conditions are different for SDP problems with different structures. For problem \eqref{SDPLR}, it is $r> \lceil \sqrt{2m} \rceil.$}. These results provide the theoretical guarantee for the BM approach to solve \eqref{SDP0}. However, when it comes to \eqref{SDP}, the above mentioned results are not applicable. This is because $\phi$ is a general convex continuously differentiable function and there are no rank conditions for \eqref{SDP} to be equivalent to (SDPR). Moreover, because $\phi$ may not be second order differentiable, we cannot use second order stationarity to derive the global optimality condition of \eqref{SDP}. In Section~\ref{Sec-optcon}, we will derive a sufficient condition, under which $(R,y)$ and $\big(\hR\hR^\top,y\circ y\big)$ are global optimal solutions of (SDPR) and (SDP), respectively. Compared to previous results, our condition does not involve the second order information of $\phi.$ This allows us to solve SDP problems using only the first order information of $\phi$. The special low rank factorization such that the first $k$ rows are fixed makes our theoretical analysis more complicated, because we have to derive an explicit formula of the dual variable of the constraint $X_{1:k,1:k}=I_k.$ Our optimality condition suggests us to solve problem \eqref{SDP} by adaptively adjusting the rank parameter $r$ and the support of $y,$ which will be elaborated in the next subsection.

\subsection{A rank-support-adaptive feasible method}\label{Rank-Adap-Alg}

In Section~\ref{Sec-Algorithm}, we will design a 
rank-support-adaptive feasible method to solve (SDP) by solving (SDPR) with an adaptive parameter $r$. 
Intuitively speaking, the rank-adaptive strategy is to dynamically increase the parameter $r$ to escape from a saddle point or
decrease it to reduce the problem size. Because (SDPR) also involves a vector variable $y,$ we will apply a support-adaptive strategy to $y,$ that is, we drop some entries of $y$ if they are nearly zero and add them back when escaping from a saddle point. The support-adaptive strategy, while hasn't been used in the literature, is useful for active set identification.

Although the rank-adaptive idea has been frequently used in the literature of 
solving low rank SDP  (see \cite{bellavia2021relaxed,Staircase1,Boumal2,Staircase2}) and matrix completion (see \cite{rankadap,lee2022escaping,rankada1,rankada2}) problems, this technique is considered as a heuristic without convergence rate analysis. In order to fill the gap between theory and practice, we will prove that our algorithm can return an approximate KKT solution under the LICQ (linear independence constraint qualification) assumption. Surprisingly, due to the special structure of (SDPR), the iteration complexity of our algorithm is better than the state-of-the-art complexity of Riemannian trust region method in \cite{BouRTR} by Boumal et al. Apart from the global sub-linear convergence, we also prove the local linear convergence of our algorithm under a quadratic growth condition, which will be shown to be shown to be a generic property for linear SDP problems with generic data $(b,C)$ in Section~\ref{Sec-randper}. In \cite{erdogdu2022convergence}, Drdogdu et al. also prove that the quadratic growth property is generic for a linear SDP with the special affine constraint $\dd(X)=e.$ Here we substantially extended their result from the special affine
constraint to the general case  $\A(X)=b.$ Our local linear convergence result requires us to identify the correct rank and support at the optimal solution $(X^*,x^*).$ This additionally implies the importance of rank-adaptive technique because it not only reduces the problem dimension but also improve the algorithm's efficiency.

\subsection{Generic smoothness and random perturbation}\label{Rand-Sing} 

The generic\footnote{In this paper, by ``generic", we always means all data points in some linear space except for a set of measure 0.} smoothness of linear SDP problem was first established in \cite{genesmo} by Alizadeh et al. by assuming $\A,b,C$ are generic. (We said that an SDP problem has the smoothness property if every feasible point satisfy the LICQ.)
In \cite{CIFrank}, Cifuentes further proved that  if the linear mapping $\A$ is generic, then any feasible point of the low rank SDP will satisfy LICQ. However, these two results are mainly of theoretical interest because in practice, the coefficient matrices $A_i$'s usually have sparse plus low rank property, i.e., they are not generic. A simple idea is to add a random small perturbation to $A_i$. However, this will destroy its sparsity and low rank property, which is vital for computational efficiency. Unlike the linear mapping $\A,$ random perturbation to the vector $b$ would not increase the computational complexity. Moreover, in \cite{Sardap}, Scholtes and St\"ohr proved that for a generic $b,$ every feasible point of a standard nonlinear optimization problem satisfies the LICQ property (see corollary 2 of \cite{Sardap}). Therefore, a random perturbation to $b$ is enough to eliminate all the singular points in (SDPR). In Section~\ref{Sec-randper}, we will use a manifold version of the Morse-Sard's Theorem \cite{DT} to improve the generic smoothness result in \cite{Sardap}. 
We  show that if the constraints of (SDPR) can be separated into two parts such that one of them constructs a manifold, then a random perturbation to the entries of $b$ in the other part will suffice to eliminate all the singular points in (SDPR) with probability one. Note that the random perturbation strategy alone cannot completely overcome the ill-conditioning issue in an  SDP problem with degenerate optimal solutions. We will further develop an adaptive preconditioning technique to make the problem 
better conditioned.


\subsection{Refinement of dual variable}
The optimality of a given point $(X,x)$ in (SDP) cannot be verified without the dual variable. In Section~\ref{Sec-RRdual}, we will discuss how to recover the dual variable from the primal solution of (SDPR). Previous works mostly focused on simple smooth SDP problems, where the dual variable can easily be recovered by solving a well-conditioned linear system (see (13) of \cite{Boumal2}). However, when the Jacobian of the constraint matrix is ill-conditioned or even singular, the dual variable computed from the above mentioned linear system may not be accurate. In order to overcome this issue, we will design a method based on ALM to improve the feasibility of the dual variable in Section~\ref{Sec-RRdual}.

\subsection{Organization of this paper}

The paper is organised as follows. In Subsection~\ref{subsec-notation}, we will present some frequently used notations. 
In Section~\ref{Sec-optcon}, we derive the optimality condition of (SDPR). In Section~\ref{Sec-Algorithm}, we present our rank-support-adaptive feasible method together with its convergence analysis. In Section~\ref{Sec-randper}, we  prove the generic smoothness of (SDPR) and describe the random perturbation and adaptive-preconditioning strategies. We also discuss the computational complexity of our feasible method. In Section~\ref{Sec-RRdual}, we develop techniques to improve the feasibility of the dual variable. In Section~\ref{Sec-numerexp}, we present 
numerical experiments to evaluate the performance of our method. We  end the paper with a brief conclusion in Section~\ref{Sec-conc}. 

\subsection{Notations}\label{subsec-notation}
Now we state some useful notations. Throughout this paper, without special indication, $n,m,r,k,p,b,R,y,\phi,\K,\A,B$ are always defined as in (SDPR). Besides, we always let $C:=\nabla_X \phi\(X,x\)$ and $c:=\nabla_x \phi\(X,x\).$ Let $\|.\|$ be the Frobenius norm and $\<A,B\>:=\Tr(AB^\top)$ be matrix inner product. Let $e$ be the column vector whose entries are all ones.
For any linear subspace $V\subset \R^n,$ $V^\perp$ is the orthogonal complement of $V$ in $\R^n.$ For any matrices $A,B$ with the same number of columns, $[A;B]$ denotes the 
matrix obtained by appending $B$ to the last row of $A$. For any vector $x$ in some metric space $\E$ and $\delta>0,$ ${\rm B}_\delta(x):=\left\{y\in \E:\ {\rm dist}(y,x)< \delta\right\}.$

In order to explicitly deal with the first $k\times k$ fixed block in (SDP) and the special factorization in (SDPR), we will introduce some notations and operators that will be frequently used throughout this paper. Let $C=[C_{11},C_{12};C_{21},C_{22}],$ where $C_{11}\in \S^k$ and $C_{22}\in \S^{n-k}.$ 
Denote $J_{n-k,n}=[0_{(n-k)\times k},I_{n-k}]\in\R^{(n-k)\times n},$ 
$I_{k,r} = [I_k,0_{k\times (r-k)}]\in\R^{k\times r}$ 
and $I_{n, k}=[I_k;0_{(n-k)\times k}]\in\R^{n\times k}.$
Define the linear mapping $\L_R:\R^{(n-k)\times r}\rightarrow \S^n$ such that for any $H\in \R^{(n-k)\times r}$,
\begin{equation}\label{Aug_7_2}
\L_R(H):=\begin{pmatrix}0_{k\times k}&I_{k,r}H^\top \\[3pt]
 HI_{k,r}^\top & RH^\top+HR^\top\end{pmatrix}.
\end{equation}
The above linear operator is the differential of the factorization $R\rightarrow [I_{k,r};R][I_{k,r}^\top,R^\top].$ For any matrix $M\in \R^{(n-k)\times r},$
define $\widehat{M}:=[I_{k,r};M]\in \R^{n\times r}.$ 
The adjoint mapping of $\L_R$ is given by
\begin{equation}\label{Aug_7_2.5}
\L_R^*\(S\)=2J_{n-k,n}S\hR \quad \forall \; S\in \S^n.
\end{equation}
For any linear mapping $\E:\S^n\rightarrow \R^m,$ define the linear mapping $\tilde{\E}:\S^{n-k}\rightarrow \R^m$ such that for any $X\in \S^{n-k}$, $\tilde{\E}(X):=\E\big( J_{n-k,n}^\top XJ_{n-k,n} \big).$ For any $H\in \R^{n\times r},$ we will frequently use the following result:
For $G=R+H$,
\begin{equation}\label{equ-fre}
\E\big( \widehat{G}\widehat{G}^\top \big)=\E\big(\hR\hR^\top\big)+\E\(\L_R\(H\)\)
+\tilde{\E}\(HH^\top\).
\end{equation}
 Define $\G:\R^{(n-k)\times r}\times \R^p\rightarrow \R^{m}$ such that for any $\(R,y\)\in \R^{(n-k)\times r}\times \R^p,$
\begin{equation}\label{Aug_7_0.5}
\G\(R,y\):=\K\big(\hR\hR^\top,y\circ y\big)-b.
\end{equation} 
Then the feasible set of (SDPR) with rank parameter $r$ is the following set: 
\begin{equation}
\BA:=\left\{ \(R,y\)\in \R^{(n-k)\times r}\times \R^p:\ \G\(R,y\) = 0\right\},\notag
\end{equation}
For any $\(R,y\)\in \BA,$ the Frech\'et differential mapping of $\G$ at $(R,y)$ is given by
\begin{equation}\label{Aug_7_1.5}
\D\G\(R,y\)[\(H,h\)] = \K\( \L_R(H),2y\circ h \) \quad \forall\; (H,h)\in \R^{(n-k)\times r}\times \R^p,
\end{equation}
and its adjoint mapping is 
\begin{equation}\label{Aug_7_1}
\D \G\(R,y\)^*[\lambda] = 2\big(J_{n-k,n}\A^*\(\lambda\)\hR,\;\dd\(B^\top \lambda\) y\big)
\quad \forall\; \lambda \in \R^m.
\end{equation}
We say that $\(R,y\)\in \BA$ is a regular (or smooth) point if the LICQ holds, i.e., 
$\D \G\(R,y\)^*$
is injective, or equivalently $\D\G(R,y)$ is surjective. Otherwise, we say that $\(R,y\)$ is a non-regular (or singular) point.

For any $r\in \mathbb{N}_+,$ we define the function  $f_r(R,y): \R^{(n-k)\times r}\times \R^p$ such that  $f_r(R,y)=\phi\big(\hR\hR^\top,y\circ y\big).$

\section{Optimality of (SDPR)}\label{Sec-optcon}

In this section, we derive a sufficient optimality condition for $\big(\hR\hR^\top,y\circ y\big)$ to be an optimal solution of (SDP), where $(R,y)$ is a feasible point of (SDPR). Our optimality condition is presented in the following lemma.
\begin{lemma}\label{Optlem}
Suppose $(R,y)\in \R^{(n-k)\times r}\times \R^p$ and there exists $\lambda\in \R^m$ that satisfies the following conditions
\begin{itemize}
\item[(i)]  $\K\big(\hR\hR^\top,y\circ y\big)=b.$
\item[(ii)] $J_{n-k,n}\( C-\A^*(\lambda)\)\hR=0,$ $\(c-B^\top\lambda\)\circ y=0.$
\item[(iii)]   $C_{22}-\tilde{\A}^*(\lambda)\in \S^{n-k}_+,$ \; $c-B^\top \lambda\in \R^p_+.$
\end{itemize}
Then $\big(\hR\hR^\top,y\circ y,\lambda,\Lambda\big)$ is a KKT solution of (SDP), where 
\begin{equation}\label{Aug_8_6.5}
\Lambda:=C_{11}+I_{k,r}R^\top J_{n-k,n}\(C-\A^*(\lambda)\)I_{n,k}-I_{n,k}^\top \A^*(\lambda) I_{n,k}
\end{equation}
is the dual variable of the constraint $X_{1:k,1:k}=I_k.$ Thus, $(R,y)$ is a global minimizer of (SDPR).
\end{lemma}

\begin{proof}
From (i), we have that $(\hR\hR^\top,y\circ y)$ satisfies the primal feasibility of (SDP). Let $Q= C_{22} - \tA^*(\lambda).$
Define ${\bf S}\in \S^n$ as  follows:
\begin{equation}\label{Aug_8_8}
{\bf S}:=\begin{pmatrix} -I_{k,r}R^\top \\ I_{n-k}\end{pmatrix} Q\begin{pmatrix} -RI_{k,r}^\top & I_{n-k}\end{pmatrix}.
\end{equation}
From (iii), we have that $Q\in \S^{n-k}_+$, and hence ${\bf S}\in \S^n_+.$ After some algebraic manipulations, we get
\begin{equation}\label{Aug_8_10}
{\bf S}:=\begin{pmatrix} I_{k,r}R^\top Q RI_{k,r}^\top & -I_{k,r}R^\top Q
\\ -QRI_{k,r}^\top & Q
 \end{pmatrix}.
\end{equation}
Expanding the first equation in (ii), we get the following equation
\begin{equation}\label{Aug_8_11}
QR=-J_{n-k,n}\(C-\A^*(\lambda)\)\begin{pmatrix} I_{k,r}\\ 0_{(n-k)\times r}\end{pmatrix}.
\end{equation}
From \eqref{Aug_8_11}, we have that 
$QRI_{k,r}^\top I_{k,r}=QR.$ This implies the following equation:
\begin{equation}\label{Aug_8_13}
Q\begin{pmatrix} -RI_{k,r}^\top & I_{n-k} \end{pmatrix}\begin{pmatrix}I_{k,r}\\ R\end{pmatrix}=0.
\end{equation}
Combining \eqref{Aug_8_8} and \eqref{Aug_8_13}, we have that ${\bf S}\hR=0.$ Moreover, from \eqref{Aug_8_11}, we get
\begin{equation}\label{Aug_8_15}
-QRI_{k,r}^\top=J_{n-k,n}\(C-\A^*(\lambda)\)I_{n,k}.
\end{equation}
Substitute \eqref{Aug_8_15} into \eqref{Aug_8_10}, we get
\begin{equation}\label{Aug_8_16}
{\bf S}:=\begin{pmatrix}
-I_{k,r}R^\top J_{n-k,n}\(C-\A^*(\lambda)\)I_{n,k} & I_{n,k}^\top\(C-\A^*(\lambda)\)J_{n-k,n}^\top \\
J_{n-k,n}\(C-\A^*(\lambda)\)I_{n,k} & C_{22}-\tA^*(\lambda)
\end{pmatrix}.
\end{equation}
From \eqref{Aug_8_16}, ${\bf S}$ can be rewritten as follows:
\begin{equation}\label{Aug_8_17}
{\bf S}:=C-\A^*(\lambda)
-I_{n,k}\,\Lambda\, I_{n,k}^\top.
\end{equation}
where 
\begin{equation} \label{Aug_8_18}
\Lambda:=C_{11}-I_{n,k}^\top \A^*(\lambda) I_{n,k}+I_{k,r}R^\top J_{n-k,n}\(C-\A^*(\lambda)\)I_{n,k}.
\end{equation}
Now, consider $\lambda$ as the dual variable of the constraint $\K(X,x)=b$ and
$\Lambda$
as the dual variable of the constraint $X_{1:k,1:k}=I_k$ in (SDP). Define ${\bf s}:=c-B^\top \lambda.$ Then $\({\bf S},{\bf s}\)$ can be viewed as the dual variable of (SDP). From (i), the second equation of (ii), the second equation of (iii), ${\bf S}\in \S^n_+$ and ${\bf S}\hR=0$, we know that the KKT conditions of (SDP) hold. Therefore, $\big(\hR\hR^\top,y\circ y,\lambda,\Lambda\big)$ is a KKT solution of (SDP). 
\end{proof}

\begin{remark}\label{remKKT}
In the proof of Lemma~\ref{Optlem}, the main task is to derive the dual variable of the constraint $X_{1:k,1:k}=I_k$, i.e., equation \eqref{Aug_8_18}. Without this constraint, conditions (i)-(iii) are exactly the KKT condition of (SDP) so there is nothing to prove. Verifying the optimality condition in Lemma~\ref{Optlem} requires us to know $\lambda$ in advance. Under the LICQ condition, the linear system (ii) of $\lambda$ has at most one solution because the mapping $\D\G(R,y)^*$ in \eqref{Aug_7_1} is injective. Thus, given $(R,y)$,
we can recover the dual variable $\lambda$ by computing a least square solution of (ii). We call $(R,y)$ a saddle point if the condition (iii) is violated. 
\end{remark}

\section{Algorithm}\label{Sec-Algorithm}
In this section, we will design a rank-support-adaptive feasible method with convergence guarantee. Because our algorithm is based on Riemannian optimization, throughout this section, we always assume that the following LICQ assumption hold. In the next section, we will discuss how to handle the singular points in $\BA$ by
a random perturbation technique.

\begin{assumption}\label{asslicq}
There exists an positive integer $\bar{r}\geq k$ such that for any  $r\geq \bar{r},$ $\BA\neq \emptyset$ and the LICQ holds at every point in $\BA.$
\end{assumption}

Before we design our algorithm, we need some preparations in the next subsection.
\subsection{Preliminaries for algorithmic design}\label{Sec-geo}
We first study the geometric structure of (SDPR). We will focus on the case where $\BA$ is a Riemannian submanifold embedded in $\R^{(n-k)\times r}\times \R^p$ (see Section 3 and 5 in \cite{Intromani} for properties of embedded manifold). 
For any $(R,y)\in \BA$, its tangent space is
\begin{equation}\label{Tansp}
{\rm T}_{(R,y)}\BA:=\left\{ (H,h)\in \R^{(n-k)\times r}\times \R^p:\ \D \G(R,y)[(H,h)]=0 \right\}
\end{equation}
where $\D\G(R,y)$ is defined in \eqref{Aug_7_1.5}. 
The projection $\P_{(R,y)}: \R^{(n-k)\times r}\times \R^p\rightarrow {\rm T}_{(R,y)}\BA$ is the orthogonal projection onto the tangent space, which is given by:
\begin{equation}\label{Proj}
\P_{(R,y)}\((D,d)\) = (D,d) - \D \G(R,y)^* \big(\D \G(R,y)\D \G(R,y)^*\big)^{-1}
 \D \G(R,y)[D,d].
 \end{equation}
 The Riemannian gradient is the projection of the Euclidean gradient onto the tangent space:
\begin{equation}\label{Riegrad}
\g f_r(R,y):=\P_{(R,y)}\( \nabla_{R,y} f_r(R,y)\)
= \nabla_{R,y} f_r(R,y)-\D \G(R,y)^*[\lambda],
\end{equation}
where $\lambda$ is the solution of the following linear system:
\begin{equation}\label{projeq}
\D \G(R,y)\D \G(R,y)^*[\lambda]  =
\D \G(R,y)[ \nabla_{R,y} f_r(R,y)].
\end{equation}
Because $(R,y)$ satisfies LICQ, $\D \G(R,y)^*$ is injective. Thus \eqref{projeq} has a unique solution. Recall that a retraction mapping $\Re_{(R,y)}: {\rm T}_{(R,y)}\BA\rightarrow \BA$ satisfies that for any $(H,h)\in {\rm T}_{(R,y)}\BA$ and $t\in \R,$
\begin{equation}\label{Retrac}
\Re_{(R,y)}\(t(H,h)\)=(R,y)+t(H,h)+o(t).
\end{equation}
Note that a retraction mapping is not unique, it can be any mappings that satisfies \eqref{Retrac}. In our algorithm, we use the Newton retraction proposed by Zhang in \cite{newretrac}, which has been proven to be a second-order retraction. The main cost in each iteration of 
the Newton retraction is in solving several linear systems of the form \eqref{projeq}.

For the purpose of convergence analysis, we need the following definition of 
an approximate stationary point.
\begin{definition}\label{appsta}
Suppose Assumption~\ref{asslicq} holds and $r\geq \bar{r}.$ For any $(R,y)\in \BA$ and $\epsilon_g,\epsilon_h>0,$ we call $(R,y)$ an $(\epsilon_g,\epsilon_h)-$stationary point of (SDPR) if 
\begin{align}
& \left\|\g f_r(R,y)\right\| \leq \epsilon_g \label{appsta1}\\
&  C_{22}-\tilde{\A}^*(\lambda)\succeq -\epsilon_h I,\
c-B^\top \lambda \geq -\epsilon_h e,
 \label{appsta2}
\end{align}
where $\lambda\in \R^m$ satisfies \eqref{projeq} and $e$ is the vector of all ones.
\end{definition}
 If $\epsilon_g=\epsilon_h=0,$ then \eqref{appsta1} and \eqref{appsta2} implies the condition (ii) and (iii) in Lemma~\ref{Optlem}. From Lemma~\ref{Optlem}, we know that $\big(\hR\hR^\top,y\circ y\big)$ is an optimal solution of (SDP). This is the motivation of Definition~\ref{appsta}. In \cite{BouRTR}, Boumal et al. use the norm of the Riemannian gradient and the positive semi-definiteness of the Riemannian Hessian to measure the second order stationarity. Our Definition~\ref{appsta} is as strong as their definition of stationarity for a linear SDP problem. This is because condition \eqref{appsta1} also uses the norm of the Riemannian gradient. Moreover, for a linear SDP problem, condition \eqref{appsta2} implies that the smallest eigenvalue of the Riemannian Hessian is lower bounded by $-2\epsilon_h$ (see (17) of \cite{Boumal2}). Therefore, it is reasonable for us perform complexity analysis using Definition~\ref{appsta}. The advantage of Definition~\ref{appsta} is that we don't need the second order information of $\phi$, so we can just assume that $\phi$ is continuously differentiable.

\subsection{Riemannian gradient descent}\label{Sec-grad}

In this subsection, we will apply the Riemannian gradient descent method 
to decrease the function value in (SDPR). We make the following assumption.
\vskip-4mm
\begin{assumption}\label{assgrad}
There exists $L_g\geq 0$ and $\rho_g>0$ such that for any $\bar{r}\leq r\leq n$ and any $(R,y)\in \BA$
\begin{multline}\label{sufides}
\left|f_r\(\Re_{(R,y)}\((H,h)\)\)-\left[ f_r(R,y)+\<(H,h),\g f_r(R,y)\>\right]\right|\\
\leq \frac{L_g}{2}\(\|H\|^2+\|h\|^2\),
\end{multline}
for all $(H,h)\in {\rm T}_{(R,y)}\BA$ such that $\sqrt{\|H\|^2+\|h\|^2}\leq \rho_g.$
\end{assumption}
Assumption~\ref{assgrad} is equivalent to Assumption 2.6 in \cite{BouRTR}, where it is used to prove the convergence of a Riemannian gradient descent method. Moreover, from Lemma 2.7 in \cite{BouRTR}, Assumption~\ref{assgrad} holds when $\BA$ is a compact submanifold of $\R^{(n-k)\times r}\times \R^p,$ the retraction mapping is globally defined and $f_r$ has Lipschitz continuous gradient in the convex hull of $\BA.$ With Assumption~\ref{assgrad}, we have the following proposition about the sufficient decrease of function value. We omit the proof because it is the same as Theorem 2.8 in \cite{BouRTR}.

\begin{proposition}\label{propgrad}
Suppose Assumption~\ref{asslicq} and Assumption~\ref{assgrad} hold. Then there exists $\gamma_g>0$ such that for any $\bar{r}\leq r\leq n$ and $(R,y)\in \BA,$ there exists $t>0$ satisfying the following inequality: 
\vskip-4mm
\begin{multline}\label{Sep_25_1}
f_r\(\Re_{(R,y)}\(-t\cdot\g f_r(R,y)\)\)-f_r(R,y)\\
\leq -\gamma_g\| \g f_r(R,y)\|\min\{1,\| \g f_r(R,y) \|\}.
\end{multline}
\end{proposition}

Note that the parameter $t$ in Proposition~\ref{propgrad} is unknown in advance. However, from Lemma 10 in \cite{BouRTR}, we can find such a stepsize by applying Armijo line-search. 

\subsection{Escaping from a saddle point}\label{Sec-hess}
Proposition~\ref{propgrad} shows that the function value can be decreased sufficiently when the optimality condition \eqref{appsta1} is violated. In this subsection, we show that when the second optimality condition \eqref{appsta2} is violated, we can also decrease the function value sufficiently by updating the rank and support of $R$ and $y$ respectively. We start by defining the following set:
\begin{equation}\label{newset}
\BA^\tau:=\left\{ (S,Y)\in \R^{(n-k)\times (r+\tau)}\times \R^{p\times 2}:\ \K\big( \hS\hS^\top,\dd\(YY^\top\) \big)-b =0\right\}.
\end{equation}
The above set $\BA^\tau$ is constructed by adding $\tau$ columns to the matrix variable and 1 column to the vector variable of $\BA$. These new columns are for us to increase the rank and support. The main idea of the support updating strategy is to consider every entry of $y$ as an $1\times 1$ matrix block and apply the rank adaptive strategy to them independently. This is why we add a new column to the vector variable. Note that we can maintain the vector structure by introducing $y:=\sqrt{\dd\(YY^\top\)}$ after updating the support. Now we are ready to show how to decrease the function value by adding new columns. 

We first study the LICQ property of $\BA^\tau.$ The following definitions are similar to those given in \eqref{Aug_7_0.5}, \eqref{Aug_7_1.5} and \eqref{Aug_7_1}. Define $\G_\tau:\R^{(n-k)\times (r+\tau)}\times \R^{p\times 2}\rightarrow \R^m$ such that for any $(S,Y)\in \R^{(n-k)\times (r+\tau)}\times \R^{p\times 2},$
\begin{equation}\label{Sep_25_7}
\G_\tau(S,Y):=\K\big( \hS\hS^\top,\dd\(YY^\top\) \big)-b.
\end{equation}
We have that for any $\(\lambda,U,V\)\in \R^m\times \R^{(n-k)\times (r+\tau)}\times \R^{p\times 2},$
\vskip-4mm
\begin{equation}\label{Sep_25_9}
\D \G_\tau(S,Y)[(U,V)]=\K\( \L_S(U),\;2\dd(YV^\top) \),
\end{equation}
and its adjoint mapping is given by
\begin{equation}\label{Sep_25_8}
\D \G_\tau(S,Y)^*[\lambda]=2\big( J_{n-k,n} \A^*(\lambda)\hS,\;\dd\(B^\top \lambda\)Y\big).
\end{equation}
For any $(S,Y)\in \BA^\tau,$ $\big(S,\sqrt{\dd(YY^\top)}\big)\in \Omega_{r+\tau}.$  Moreover, from \eqref{Aug_7_1} and \eqref{Sep_25_8}, $\D \G_\tau(S,Y)^*$ is an injective mapping if and only if $\D \G \big(S,\sqrt{\dd\(YY^\top\)}\big)^*$ is an injective mapping. Therefore, from Assumption~\ref{asslicq}, we know that for any $r\in \mathbb{N}_+$ such that $r\geq \bar{r},$ every point $(S,Y)\in \BA^\tau$ satisfies the LICQ property. 

With the LICQ property, we move on to study the geometric property of $\BA^\tau.$ The tangent space of $\BA^\tau$ at $(S,Y)$ is given as follows:
\begin{equation}\label{Sep_25_10}
{\rm T}_{(S,Y)}\BA^\tau:=\left\{ (U,V)\in \R^{(n-k)\times (r+\tau)}\times \R^{p\times 2}:\ \D \G_\tau(S,Y)[(U,V)]=0\right\}.
\end{equation}
Because of the LICQ property,  $\BA^\tau$ is a Riemannian submanifold of $\R^{(n-k)\times (r+\tau)}\times \R^{p\times 2}.$ There exists a second order retraction $\Re_{(S,Y)}:{\rm T}_{(S,Y)}\BA^\tau\rightarrow \BA^\tau$ such that for any $(U,V)\in {\rm T}_{(S,Y)}\BA^\tau$ and $t\in \R,$
\vskip-4mm
\begin{equation}\label{Sep_25_11}
 \Re_{(S,Y)}\( t(U,V) \)=(S,Y)+t(U,V)+\frac{t^2}{2}(W,Z)+o(t^2),
\end{equation}
for some $(W,Z)\in \big({\rm T}_{(S,Y)}\BA^{\tau}\big)^\perp.$ 
In practice, one can use the metric projection (see 5.48 of \cite{Intromani}), the Newton retraction in \cite{newretrac} or the orthographic retraction in \cite{Pretrac} as a second order retraction. Since $(W,Z)\in \big({\rm T}_{(S,Y)}\BA^{\tau}\big)^\perp,$ there exists $\hat{\lambda}\in \R^m$ such that 
\begin{equation}\label{Sep_25_12}
(W,Z)=\D \G_\tau(S,Y)^*[\hat{\lambda}]=2\big( J_{n-k,n}\A^*(\hat{\lambda})\hS,\;\dd\big(B^\top\hat{\lambda}\big)Y \big).
\end{equation}
Substituting \eqref{Sep_25_11} and \eqref{Sep_25_12} into the equation $\G_\tau\(  \Re_{(S,Y)}\( t(U,V) \) \)=0$, we obtain the following equation, which says that the coefficient of $t^2$ in the expansion of $\G_\tau\(  \Re_{(S,Y)}\( t(U,V) \) \)$ is zero:
\begin{equation}\label{Sep_25_13}
\frac{1}{2}\D\G_\tau(S,Y)\left[ \D \G_\tau(S,Y)^*[\hat{\lambda}]\right]+\tilde{\K}\(UU^\top,\dd(VV^\top)\)=0.
\end{equation}
The above equation also implies that $\hat{\lambda}$ is the unique solution of the linear system \eqref{Sep_25_13}. Note that we have used the formula \eqref{equ-fre} to derive the above equation.
Therefore, $(W,Z)$ is a smooth mapping of $(U,V).$ Define $f_{r,\tau}:\R^{(n-k)\times (r+\tau)}\times \R^{p\times 2}\rightarrow \R$ such that for any $(S,Y)\in \R^{(n-k)\times (r+\tau)}\times \R^{p\times 2}$
\begin{equation}\label{Sep_25_14}
f_{r,\tau}\((S,Y)\):=\phi\big( \hS\hS^\top,\dd\(YY^\top\) \big).
\end{equation}
Now, with the geometric structure of $\BA^\tau,$ we are able to show how the function value will change by adding new columns $(H,h)\in \R^{(n-k)\times \tau}\times \R^p$ to $(R,y)\in \BA.$ In order to do so, consider the following variables:
\begin{equation}\label{Sep_25_14.5}
(S,Y):=\([R,0_{(n-k)\times \tau}],[y,0_{p\times 1}]\),\ (U,V):=\( [0_{(n-k)\times r},H],[0_{p\times 1},h] \)
\end{equation}
 such that $H\in \R^{(n-k)\times \tau}$ and $h\in \R^p.$ It is easy to see that $(S,Y)\in \BA^\tau$ and $(U,V)\in {\rm T}_{(S,Y)}\BA^\tau.$ Also, from \eqref{Sep_25_12}, we know that the last $\tau$ columns of $W$ and the second column of $Z$ are zero columns respectively i.e., $W=[W_1,0_{(n-k)\times \tau}]$ and $Z=[z_1,0_{p\times 1}]$ for some $W_1\in \R^{(n-k)\times r}$ and $z_1\in \R^p.$ 
Because $\phi$ is continuously differentiable, $f_{r,\tau}$ is also continuously differentiable. From \eqref{Sep_25_11}, we have the following result:
\vskip-4mm
\begin{align}
&f_{r,\tau}\( \Re_{(S,Y)}\( t(U,V) \) \)=f_{r,\tau}\Big((S,Y)+t(U,V)+\frac{t^2}{2}(W,Z)\Big)+o(t^2)\notag \\
&=f_{r,\tau}\(\Big(\big[ R+\frac{t^2}{2}W_1,tH \big],\big[y+\frac{t^2}{2}z_1,th\big]\Big)\)+o(t^2)\notag \\
&=\phi\Big( \hR\hR^\top+\frac{t^2}{2}\L_R(W_1)+t^2J_{n,n-k}^\top HH^\top J_{n,n-k},\ y\circ y+t^2 y\circ z_1+t^2h\circ h \Big)+o(t^2).\notag \\
&= f_r\(R,y\)+\frac{t^2}{2}\big(\< C,\L_R(W_1) \>+\<c,2y\circ z_1\>\big)+t^2\(\<C_{22},HH^\top\>+\<c,h\circ h\>\)+o(t^2),\label{Sep_25_15}
\end{align}
where $C = \nabla_X\phi(RR^\top,y\circ y)$ and $c=\nabla_x\phi(RR^\top,y\circ y)$.
By using \eqref{Aug_7_2.5}, we have that
\begin{align}
&\< C,\L_R(W_1) \>+\<c,2y\circ z_1\>= \<2J_{n-k,n}C\hR,W_1\>+\<2c\circ y,z_1\> \notag \\
&=\< \nabla_{R,y} f_r\(R,y\),\(W_1,z_1\) \>=\< \g f_r\( R,y \)+\D \G(R,y)^*[\lambda],\(W_1,z_1\) \>\notag \\
&=\< \g f_r\( R,y \),\(W_1,z_1\) \>+\< \lambda,\D \G(R,y)[\(W_1,z_1\)] \>,\label{Sep_25_16}
\end{align} 
where the third equality comes from \eqref{Riegrad} and $\lambda$ is the unique solution of \eqref{projeq}. From \eqref{Sep_25_12} and the relation between $(W,Z)$ and $(W_1,z_1)$, we know that $\(W_1,z_1\)=\D \G(R,y)^*[\hat{\lambda}]\in \big({\rm T}_{(R,y)}\BA\big)^\perp.$ Since $\g f\( R,y \)\in {\rm T}_{(R,y)}\BA,$ we have 
\begin{equation}\label{Sep_25_17}
\< \g f_r\( R,y \),\(W_1,z_1\) \>=0.
\end{equation}
 Also, from \eqref{Sep_25_12} and \eqref{Sep_25_13}, we know that 
\begin{align}
&\D \G(R,y)[\(W_1,z_1\)]=\D \G(R,y)\left[\D \G(R,y)^*[\hat{\lambda}]\right]=\D \G_\tau(S,Y)\left[\D \G_\tau(S,Y)^*[\hat{\lambda}]\right]\notag \\
&=-2\tilde{\K}\(UU^\top,\dd(VV^\top)\)=-2\tilde{\K}\(HH^\top,h\circ h\),\label{Sep_25_18}
\end{align}
where the first equality comes from $\(W_1,z_1\)=\D \G(R,y)^*[\hat{\lambda}]$, and the second equality comes from $\(S,Y\)=\([R,0_{(n-k)\times \tau}],[y,0_{p\times 1}]\).$
Substituting \eqref{Sep_25_17} and \eqref{Sep_25_18} into \eqref{Sep_25_16}, we get
\begin{equation}\label{Sep_25_19}
\< C,\L_R(W_1) \>+\<c,2y\circ z_1\>=-2\<\lambda,\tilde{\K}\(HH^\top,h\circ h\)\>.
\end{equation}
Substituting \eqref{Sep_25_19} into \eqref{Sep_25_15}, we get
\begin{align}
&f_{r,\tau}\( \Re_{(S,Y)}\( t(U,V) \) \)\notag \\
&=f_r(R,y)+t^2\(\<C_{22},HH^\top\>+\<c,h\circ h\>-\<\lambda,\tilde{\K}\(HH^\top,h\circ h\)\>\)+o(t^2)\notag \\
&=f_r(R,y)+t^2\(\< C_{22}-\tilde{\A}^*(\lambda),HH^\top \>+\<c-B^\top \lambda,h\circ h\>\)+o(t^2). \label{Sep_25_20}
\end{align}
Note that \eqref{Sep_25_20} is a local property at $(R,y)\in \BA.$ Similar to Assumption~\ref{assgrad}, we make the assumption below, which says that \eqref{Sep_25_20} can be extended globally on $\BA.$

\begin{assumption}\label{asshess}
Consider $\tau\in \mathbb{N}_+.$ There exists $\rho_h>0$ and $L_h\geq 0$ such that for any 
integer
$r$ satisfying $\bar{r}\leq r\leq n,$ $(R,y)\in \BA,$ 
\begin{align*}
&\left|f_{r,\tau}\( \Re_{(S,Y)}\( (U,V) \) \)- \(f_r(R,y)+\< C_{22}-\tilde{\A}^*(\lambda),HH^\top \>+\<c-B^\top \lambda,h\circ h\>\)\right|\notag \\
&\leq \frac{L_h\(\|H\|^2+\|h\|^2\)^{\frac{3}{2}}}{6},\label{Sep_26_1}
\end{align*}
for all $(H,h)\in \R^{(n-k)\times \tau}\times \R^p$ such that $\sqrt{\|H\|^2+\|h\|^2}\leq \rho_h,$ where $(S,Y)$ and $(U,V)$ are defined as in \eqref{Sep_25_14.5} and $\lambda$ is the unique solution of \eqref{projeq}.
\end{assumption}

Assumption~\ref{asshess} is similar to Assumption 3.2 in \cite{BouRTR}. Unlike \cite{BouRTR}, we avoid using the Riemannian Hessian. Also, we fix our search space to be $\R^{(n-k)\times \tau}\times \R^p$ instead of the whole tangent space ${\rm T}_{(S,Y)}\BA^\tau.$ We remark that Assumption~\ref{asshess} will hold when $\phi$ has locally Lipschitz gradient and $\BA^\tau$ is a compact smooth submanifold with a globally defined second order retraction i.e., one can replace $o(t^2)$ with $c_h t^3$ for some constant $c_h>0$ in \eqref{Sep_25_11}, when $\sqrt{\|U\|^2+\|V\|^2}=1$ and $t\in [-\rho_h,\rho_h].$ With Assumption~\ref{asshess}, we can prove the following proposition on how to escape from a saddle point.

\begin{proposition}\label{prophess}
Suppose Assumption~\ref{asslicq} and Assumption~\ref{asshess} hold. Then there exists $\gamma_h>0$ such that for any integer $r$ satisfying $\bar{r}\leq r\leq n$ and $(R,y)\in \BA,$ there exists $t\geq 0$ satisfying the following inequality:
\begin{equation}\label{Sep_26_2}
f_{r,\tau}\(\Re_{(S,Y)}\( t(U,V) \) \)-f_r(R,y)\leq \gamma_h\xi\min\{1,\xi^2\},
\end{equation}
where $(S,Y)=\([R,0_{(n-k)\times \tau}],[y,0_{p\times 1}]\)$, $(U,V)=\( [0_{(n-k)\times r},H],[0_{p\times 1},h]\)$ 
such that $(H,h)\in \R^{(n-k)\times \tau}\times \R^p$ is an optimal solution of the following problem:
\begin{equation}\label{escprob}
\xi:=\min\left\{ \<C_{22}-\tilde{\A}^*(\lambda),HH^\top\>+\<c-B^\top\lambda,h\circ h\>:\ \|H\|^2+\|h\|^2=1 \right\},
\end{equation}
where $\lambda$ is the unique solution of \eqref{projeq}.
\end{proposition}

\begin{proof}
If $\xi\geq 0,$ then let $t=0$ and \eqref{Sep_26_2} holds. Suppose $\xi<0.$ From the inequality in Assumption \ref{asshess}
 and \eqref{escprob}, we have that
\begin{equation}\label{Sep_26_3}
f_{r,\tau}\(\Re_{(S,Y)}\( t(U,V) \) \)-f_r(R,y)\leq \xi t^2+\frac{L_h}{6}t^3,
\end{equation}
for any $t\in [-\rho_h,\rho_h].$ For any constant $t_0>0,$ $\gamma_h'\in (0,1)$ and $\delta\in (0,1),$ let $k$ 
be the smallest number in $\mathbb{N}$ such that $t_k:=t_0\delta^k$ satisfies the following inequality
\begin{equation}\label{Sep_26_4}
f_{r,\tau}\(\Re_{(S,Y)}\( t_k(U,V) \) \)-f_r(R,y)\leq \gamma_h'\xi t_k^2.
\end{equation}
The existence of $k$ is guaranteed by \eqref{Sep_26_3}. If $k=0,$ then we have that
\begin{equation}\label{Sep_26_5}
f_{r,\tau}\(\Re_{(S,Y)}\( t_k(U,V) \) \)-f_r(R,y)\leq \gamma_h'\xi t_0^2.
\end{equation}
If $k>0,$ because $k$ is the smallest integer satisfying \eqref{Sep_26_4}, $t_k/\delta$ doesn't satisfy \eqref{Sep_26_4}. This together with \eqref{Sep_26_3} implies that
\begin{equation}\label{Sep_26_6}
t_k/\delta\geq \min\left\{ \rho_h,\ {6(\gamma_h'-1)\xi}/{L_h}\right\}.
\end{equation}
Substituting \eqref{Sep_26_6} into \eqref{Sep_26_4}, we get
\begin{equation}\label{Sep_26_7}
f_{\tau,r}\(\Re_{(S,Y)}\( t_k(U,V) \) \)-f_r(R,y)\leq \gamma_h'\delta^2\min\left\{ \rho_h^2,{36 (\gamma_h'-1)^2\xi^2}/{L_h^2} \right\}\xi.
\end{equation}
Combining \eqref{Sep_26_5} and \eqref{Sep_26_7}, we have that
\begin{equation}\label{Sep_26_8}
f_{r,\tau}\(\Re_{(S,Y)}\( t_k(U,V) \) \)-f_r(R,y)\leq \gamma_h\min\{1,\xi^2\}\xi,
\end{equation}
where $\gamma_h:=
\gamma_h' \min\left\{ t_0^2,\delta^2 \rho_h^2,36 (\gamma_h'-1)^2\delta^2/L_h^2 \right\}.$
\end{proof}
Similar to Proposition~\ref{propgrad}, the proof of Proposition~\ref{prophess} suggests that we can use Armijo line-search to decrease to the function value sufficiently without knowing $L_h$ and $\rho_h$ in advance. Note that after we get $(S_1,Y_1):=\Re_{(S,Y)}\( t(U,V) \),$ we define $(R_1,y_1):=(S_1,\sqrt{\dd(Y_1Y_1^\top)}).$ Then we have that $(R_1,y_1)\in \Omega_{r+\tau}$ and $f_{r+\tau}(R_1,y_1)=f_{r,\tau}(S_1,Y_1).$
 
\subsection{A rank-support-adaptive feasible method and its convergence analysis}\label{Sec-alg}
Now, we are ready to state our rank-support-adaptive feasible method in Algorithm~\ref{alg1}.

\begin{algorithm}
\caption{A rank-support-adaptive feasible method for (SDP)}
\label{alg1}
\begin{algorithmic}
\STATE{{\bf Parameters:} $r_0\in \mathbb{N}_+\cap[\bar{r},n],$ $\tau\in \mathbb{N}_+,$ $\epsilon_g>0$, $\epsilon_h>0$, $\{\epsilon_i\}_{i\in \mathbb{N}_+}\subset \R_+$}
\STATE{{\bf Input:} $(R_0,y_0)\in \Omega_{r_0}$}
\STATE{{\bf Initialization:} $i\gets 0$}
\WHILE{{\bf the termination criterion is not met:} }
\STATE{Obtain $\lambda_i$ from \eqref{projeq}}
\IF{$\| \g f_{r_i}(R_i,y_i) \|\geq \epsilon_g$}
\STATE{Obtain $t_i\in \R_+$ satisfying \eqref{Sep_25_1} \COMMENT{Prop.~\ref{propgrad}: gradient descent}}
\STATE{$(R_{i+1},y_{i+1})\gets\Re_{(R_i,y_i)}\(-t_i\cdot\g f_{r_i}(R_i,y_i)\)$}
\STATE{$r_i^+\gets r_i$}
\ELSE
\STATE{Obtain $\xi_i,H_i,h_i$ from \eqref{escprob}}
\IF{$\max\{-\xi_i,0\}\leq \epsilon_h$}
\RETURN$\(R_i,y_i,\lambda_i\)$
\ELSE
\STATE{Obtain $S_i,Y_i,U_i,V_i$ from \eqref{Sep_25_14.5}}
\STATE{Obtain $t_i\in \R$ satisfying \eqref{Sep_26_2} \COMMENT{Prop.~\ref{prophess}: escaping saddle point}}
\STATE{$(S_{i+1},Y_{i+1})\gets \Re_{(S_i,Y_i)}\( t_i(U_i,V_i) \) $}
\STATE{$(R_{i+1},y_{i+1})\gets \big(S_{i+1},\sqrt{\dd\(Y_{i+1}Y_{i+1}\)}\,\big)$ \COMMENT{update rank and support}}
\STATE{$r_i^+\gets r_i+\tau$}
\ENDIF
\ENDIF
\IF{$R_{i+1}$ {\rm has more than} $n$ {\rm columns}}
\STATE{Obtain $R\in \R^{(n-k)\times n}$ satisfying $\hR\hR^\top=\hR_{i+1}\hR_{i+1}^\top$ \COMMENT{Control rank}}
\STATE{$R_{i+1}\gets R$}
\ENDIF
\STATE{Find $(R',y')\in \Omega_{r'}$ such that $\bar{r}\leq r'\leq r_i^+),$ ${\rm supp}(y')\subset {\rm supp}(y_{i+1})$ and $f_{r'}(R',y')\leq f_{r_i^+}(R_{i+1},y_{i+1})+\epsilon_i.$ \COMMENT{Reduce rank and support}}
\STATE{$(R_{i+1},y_{i+1})\gets (R',y')$}
\STATE{ $i\gets i+1,$ $r_i\gets r_i^+$}
\ENDWHILE

\end{algorithmic}
\end{algorithm}

Algorithm~\ref{alg1} attempts to find an $(\epsilon_g,\epsilon_h)-$stationary point of (SDPR). In the $i$th iteration, if $\| \g\( R_i,y_i \) \|$ is large, we apply the gradient descent in subsection~\ref{Sec-grad} to decrease to function value. A sufficient decrease in the 
function value is guaranteed by \eqref{Sep_25_1}. Otherwise if  $\xi_i$ defined in \eqref{escprob} is negative and small, i.e., $(R_i,y_i)$ is a nearly strict saddle point, we apply the escaping saddle point strategy in subsection~\ref{Sec-hess} to increase $r$ by $\tau$ and decrease the function value. Note that $\max\{-\xi,0\}\leq \epsilon_h$ is equivalent to the condition~\eqref{appsta2}. Similarly, \eqref{Sep_26_2} implies that we can decrease the function value sufficiently. In every iteration, we can reduce the rank and support of the iterate $(R_{i+1},y_{i+1})$ to be $(R',y')$ provided that the increase of the function value is controlled by $\epsilon_i$. Our reduction technique is motivated by \cite{rankadap}. In detail, let $\sigma_1\geq \sigma_2\geq \ldots \geq \sigma_k>0$ be the nonzero singular values of $R_{i+1}.$ Let $r':=\min\{\arg\max\{\sigma_i/\sigma_{i+1}:\ i\in [k-1]\}\}.$ If $\sigma_{r'}/\sigma_{r'+1}>\kappa_1,$ where $\kappa_1>1$ is the threshold of the rank reduction, then we drop all the singular values smaller than $\sigma_{r'}$ to get a new matrix $R_{i+1}'\in \R^{(n-k)\times r'}.$ Also, let $x:=y_{i+1}\circ y_{i+1}$ and define the set $\I:=\left\{ j\in [p]:\ x_j\leq  x_{\max}/\kappa_2 \right\},$ where $x_{\max}$ is the largest entry of $x$ and $\kappa_2>1$ is the threshold for support reduction. We let all entries of $y_{i+1}$ in $\I$ to be zero to get $y_{i+1}'.$ After this, we apply Newton retraction to $(R_{i+1}',y_{i+1}')$ to get $(R',y')\in \Omega_{r'}.$ The dropped entries of $y$ will always be zero in the Riemannian gradient descent step (see \eqref{Aug_7_1}). This zero-preserving property benefits our algorithmic implementation, because we don't have to set the small entries to be zero in every iteration. Also, 
one does not have to worry about that the rank and support are over-reduced because the escaping saddle point step allows us to increase the rank and support. 

Although the rank reduction technique works quite well in practice and the search rank $r$ is usually within a small constant factor of the rank of the optimal solution, we cannot theoretically exclude the possibility that $r$ grows to $n$
during the course of the algorithm
even if the rank of the optimal solution $r^*\ll n$. Recently, Zhang \cite{zhang2022improved} established the global guarantee of low-rank BM approach for general unconstrained nonlinear SDP problems. In detail, he proved that in the unconstrained case, for a twice differentiable, $\mu$-strongly convex and $L-$Lipschitz function $\phi$, if $r>\frac{1}{4}\(L/\mu-1\)^2r^*,$ then a second order stationary point of (SDPR) is also a global optimal solution. 
Despite the different problem settings, this result undoubtedly offers some theoretical support for our empirical observations 
that we only need to increase $r$ to a small constant factor of $r^*$ in our
algorithm.

Now we state and prove our convergence result.
\begin{theorem}\label{conv}
Suppose Assumption~\ref{asslicq}, Assumption~\ref{assgrad} and Assumption~\ref{asshess} hold. Moreover, assume that $\sum_{i=0}^\infty\epsilon_i<M$ for some constant $M>0$ and problem (SDP) is solvable with optimal value $\phi^*.$ Then for any $\epsilon_g,\epsilon_h>0,$ Algorithm~\ref{alg1} will return an $(\epsilon_g,\epsilon_h)-$stationary within $O(\epsilon_g^{-2}+\epsilon_h^{-3})$ iterations.
\end{theorem}

\begin{proof}
Suppose Algorithm~\ref{alg1} has not terminated in $K$ iterations. Let $\fval_i$ be the function value at the $i$th iteration. For any $0\leq i<K,$ there are two cases to consider. 
If $\| \g f_{r_i}(R_i,y_i) \|\geq \epsilon_g,$ then the gradient descent step will be conducted. From \eqref{Sep_25_1}, we have that  
\begin{equation}\label{Sep_26_9}
\fval_{i+1}\leq \fval_i-\gamma_g\epsilon_g\min\{1,\epsilon_g\}+\epsilon_i.
\end{equation}
If $\| \g f_{r_i}(R_i,y_i) \|< \epsilon_g.$ Then we know that $\xi_i\leq -\epsilon_h.$ From \eqref{Sep_26_2} and the step $\(R_{i+1},y_{i+1}\)\gets \big(S_{i+1},\sqrt{\dd\(Y_{i+1}Y_{i+1}^\top\)}\,\big)$ in Algorithm~\ref{alg1}, we have that
\begin{equation}\label{Sep_26_10}
\fval_{i+1}\leq \fval_i-\gamma_h\epsilon_h\min\{1,\epsilon_h^2\}+\epsilon_i.
\end{equation}
Combining \eqref{Sep_26_9} and \eqref{Sep_26_10}, we have that 
\begin{eqnarray}
\fval_K-\fval_0&=&\mbox{$\sum_{i=0}^{K-1}$} \(\fval_{i+1}-\fval_i\)
\nonumber \\
&\leq& -K\min\left\{  \gamma_g\epsilon_g\min\{1,\epsilon_g\},\gamma_h\epsilon_h\min\{1,\epsilon_h^2\} \right\}
+ \mbox{$\sum_{i=0}^{K-1}$} \epsilon_i
\nonumber \\
&\leq&
 -K\min\left\{  \gamma_g\epsilon_g\min\{1,\epsilon_g\},\gamma_h\epsilon_h\min\{1,\epsilon_h^2\} \right\}+M.
\label{Sep_26_11}
\end{eqnarray}
From \eqref{Sep_26_11} and $\fval_K\geq \phi^*,$ we have that
\begin{align}
&K\leq \frac{M+\fval_0-\fval_K}{\min\left\{  \gamma_g\epsilon_g\min\{1,\epsilon_g\},\gamma_h\epsilon_h\min\{1,\epsilon_h^2\} \right\}}\notag \\
&\leq \(M+\fval_0-\phi^*\)\( \frac{1}{\gamma_g\epsilon_g\min\{1,\epsilon_g\}}+\frac{1}{\gamma_h\epsilon_h\min\{1,\epsilon_h^2\}}\) =O(\epsilon_g^{-2}+\epsilon_h^{-3}).
\label{Sep_26_12}
\end{align}
Therefore, $K=O(\epsilon_g^{-2}+\epsilon_h^{-3}).$
\end{proof}
Our convergence rate $O\(\epsilon_g^{-2}+\epsilon_h^{-3}\)$ is better than the rate
of $O\(\epsilon_g^{-2}\epsilon_h^{-1}+\epsilon_h^{-3}\)$ in Theorem 3.9 of \cite{BouRTR}. This is because our problem (SDPR) has some special structure so that the decrease of function value in Proposition~\ref{prophess} is a function of $\xi$, which is not related to the Riemannian gradient. Note that Theorem~\ref{conv} requires that $\sum_{i=0}^\infty \epsilon_i<\infty.$ This can be achieved in many ways such as increasing the thresholds for rank-support reduction after every iteration and limit the number of rank-support reduction steps. In our inplementation, we set $\kappa_1=10$ and $\kappa_2=10000.$ The number of rank-support reductions is always smaller than 20 and so the condition is satisfied.

\subsection{Local linear convergence of Algorithm~\ref{alg1}}\label{subsec-linear}
The previous subsection establishes the global sub-linear convergence of  Algorithm~\ref{alg1}. In this subsection, we move on to prove the local linear convergence of this algorithm. We need the following quadratic growth condition, which is an extension of the rank$-r$ quadratic growth condition in \cite{bai2019proximal} by Bai et al. 
\begin{assumption}\label{rquad}
Suppose $(R^*,y^*)\in \BA$ is a $(0,0)-$stationary point of (SDPR), $(R^*,y^*)$ is said to satisfy the $\(\alpha,\epsilon\)-$quadratic growth condition if $\rr(R^*)=r$ and
\begin{equation}\label{quagrowth}
\phi\(X,x\)\geq  \phi\(X^*,x^*\)+\frac{\alpha}{2}\(\| X-X^* \|^2+\|x-x^*\|^2\),
\end{equation}
where $\(X^*,x^*\)=\big(\hR^*\hR^{* \top},y^*\circ y^*\big)$ and $\(X,x\)
=\big(\hR\hR^\top,y\circ y\big)$ for any $(R,y)\in \BA$ such that ${\rm supp}(y)={\rm supp}(y^*)$ and $\max\{ \|R-R^*\|,\|y-y^*\| \}<\epsilon.$ 
\end{assumption}

\begin{remark}\label{remquad}
We remark that for a linear SDP problem, from Theorem 3.137 in \cite{bonnans2013perturbation}, the quadratic growth condition is satisfied if and only if the second order sufficient condition is satisfied. From Proposition 15 of \cite{chan2008constraint}, the second order sufficient condition for a linear SDP
is equivalent to the dual non-degeneracy condition provided that the dual variable is unique. 
\end{remark}

With the above quadratic growth condition, we have the following lemma, which can be viewed as a manifold version of the Polyak-Lojasiewicz condition \cite{karimi2016linear,lojasiewicz1963topological}.

\begin{lemma}\label{PLlem}
Suppose Assumption~\ref{asslicq} holds, $\bar{r}\leq r\leq n$ and that $(R^*,y^*)\in \BA$ satisfies the $(\alpha,\epsilon)-$quadratic growth condition with parameters $\alpha,\epsilon>0$. Suppose further that $\g f_r$ is locally Lipschitz on $\BA.$ Then there exists $\delta,\eta>0$ such that for any $(R,y)\in \BA$ satisfying $\max\left\{ \|R-R^*\|,\|y-y^*\| \right\}<\delta$ and ${\rm supp}(y)={\rm supp}(y^*),$ the following condition holds
\begin{equation}\label{PLineq}
\|\g f_r(R,y)\|^2\geq \eta \( f_r(R,y)-f_r(R^*,y^*) \).
\end{equation}
\end{lemma}
\begin{proof}
Let $\(X^*,x^*\)=\big( \hR^*\hR^{* \top},y^*\circ y^* \big)$ and $\(X,x\)=\big( \hR\hR^{\top},y\circ y \big).$ Let $\sigma_r(R^*)$ be the smallest singular value of $R^*.$ Because $\rr(R^*)=r,$ we have that $\sigma_r(R^*)>0.$  Let $\sigma_{\min}(y^*)$ be the smallest positive value of $|y^*|$ if $\|y^*\|>0$ and $0$ if $\|y^*\|=0.$ Let $\delta\in (0,\epsilon)$ be a constant such that for any $(R,y)\in \BA$ satisfying ${\rm supp}(y^*)={\rm supp}(y)$ and $\max\left\{ \|R-R^*\|,\|y-y^*\| \right\}<\delta,$ the following properties hold.
\begin{itemize}
\item[(i)]  Let $\lambda\in \R^m$ be the multiplier of $(R,y),$ which is the unique solution of \eqref{projeq}. It holds that
\begin{multline}\label{Mar_7_8}
\<C_{22}-\tA^*(\lambda),(R-R^*)(R-R^*)^\top\>\geq -\sigma_h\|R-R^*\|^2,\\
\< c-B^\top \lambda, \(y-y^*\)\circ \(y-y^*\)\>\geq -\sigma_h\|y-y^*\|^2, \hspace{2cm}
\end{multline}
where  $\sigma_h$ is a positive number defined as follows
\begin{equation}\label{Mar_9_1}
\sigma_h:=\begin{cases} 
\alpha\min\{ \sigma_r(\hR^*)^2,\sigma_{\min}(y^*)^2 \}/8 & \sigma_{\min}(y^*)>0 \\
\alpha \sigma_r(\hR^*)^2 /8 & \sigma_{\min}(y^*)=0.
\end{cases}
\end{equation}
\item[(ii)] $\|X-X^*\|^2\geq 2(\sqrt{2}-1)\sigma_r^2(\hR^*)\| R-R^* \|^2$ and $\|x-x^*\|^2\geq \sigma_{\min}(y^*)^2\|y-y^*\|^2.$ \\[-5pt]
\item[(iii)] $f_r(R,y)-f_r(R^*,y^*)\leq 2\kappa \(\|R-R^*\|^2+\|y-y^*\|^2\)$ for some constant $\kappa>0.$
\end{itemize}

\medskip
Now we explain why the above three properties can be achieved. Let $\lambda^*\in \R^m$ be the multiplier of $(R^*,y^*).$ We know that $C_{22}-\tA^*(\lambda^*)\in \S^{n-k}_+$ and $c-B^\top\lambda^*\in \R^p_+.$ From the continuity of the symmetric positive definite linear system \eqref{projeq}, we know that $\lambda\rightarrow \lambda^*$ as $(R,y)\rightarrow (R^*,y^*).$ Thus, from the continuity of eigenvalues, property \eqref{Mar_7_8} can be guaranteed as long as $\delta>0$ is small enough. For property (ii), when $\delta>0$ is sufficiently small, $R^\top R^*$ will be positive semidefinite. Thus, from Lemma 41 in \cite{ge2017no}, we have that 
\begin{equation}\label{Mar_7_10}
\|X-X^*\|^2=\|\hR\hR^\top-\hR^*\hR^{* \top}\|^2\geq 2(\sqrt{2}-1)\sigma_r^2(\hR^*)\| R-R^* \|^2.
\end{equation}
Also, since ${\rm supp}(y)={\rm supp}(y^*),$ we have that when $y$ is sufficiently close to $y^*,$ every nonzero entry of $y$ has the same sign as that of $y^*,$ which implies the following result
\begin{equation}\label{Mar_7_11}
\|x-x^*\|^2=\| y\circ y-y^*\circ y^*\|^2 = \| (y-y^*)\circ (y+y^*) \|^2\geq \sigma_{\min}(y^*)^2\|y-y^*\|^2.
\end{equation}
For property (iii), because $\BA$ is a Riemannian manifold embedded in $\R^{(n-k)\times r}\times \R^p,$ when $\delta>0$ is sufficiently small, the geodesic $\gamma:[0,1]\rightarrow \BA$ from $(R^*,y^*)$ to $(R,y)$ satisfies that 
\begin{equation}\label{Mar_8_3.13}
\int_0^1 \| \gamma'(t) \| {\rm d}t\leq 2\sqrt{\|R-R^*\|^2+\|y-y^*\|^2}.\footnote{For a Riemannian manifold $\M$ embedded in some Euclidean space $\E$, $\lim_{x\rightarrow y}\dist_\M(y,x)/\dist_\E(y,x)=1,$ where $\dist_\M(\cdot,\cdot)$ is the geodesic distance and $\dist_\E(\cdot,\cdot)$ is the Euclidean distance. This comes from the fact that the differential of the exponential mapping at $y$ is the identity mapping (see 5.4 of \cite{manibook}).} 
\end{equation}
Moreover, from the locally Lipschitz continuity of $\g f_r$ and $\| \g f_r(R^*,y^*) \|=0,$ we have that for some constant $\kappa>0$,
\begin{equation}\label{Mar_8_3.12}
\| \g f_r(R,y) \|\leq \kappa \sqrt{\|R-R^*\|^2+\|y-y^*\|^2}.
\end{equation}
Using \eqref{Mar_8_3.13}, \eqref{Mar_8_3.12} and the Cauchy-Schwatz inequality, we have that
\begin{multline}\label{Mar_9_2}
f_r(R,y)-f_r(R^*,y^*)=\int_0^1\< \g f_r(\gamma(t)),\gamma'(t) \> {\rm d}t \\
\leq \int_0^1  \|\g f_r(\gamma(t))\|\|\gamma'(t)\| {\rm d}t\leq 2\kappa \(\|R-R^*\|^2+\|y-y^*\|^2\),
\end{multline}
which is exactly property (iii). Now, we move on to prove \eqref{PLineq}. Because $\delta\in (0,\epsilon)$ and ${\rm supp}(y)={\rm supp}(y^*),$ from Assumption~\ref{rquad}, we have the following inequality
\begin{multline}\label{Mar_7_1}
\<\nabla _X \phi(X,x),X-X^*\>+\< \nabla_x \phi(X,x),x-x^* \>\geq \phi(X,x)-\phi(X^*,x^*)\\
\geq \alpha\(\| X-X^* \|^2+\|x-x^*\|^2\)/2,
\end{multline}
where the first inequality comes from the convexity of $\phi$. Because $(X^*,x^*)$ and $(X,x)$ are feasible solutions of (SDP), we have that
\begin{equation}\label{Mar_7_2}
\<\lambda,\A(X-X^*)+B(x-x^*)\>=0,
\end{equation}
 Substituting \eqref{Mar_7_2} into the left-hand side of \eqref{Mar_7_1}, we get
\begin{multline}\label{Mar_7_3}
\<\nabla _X \phi(X,x),X-X^*\>+\< \nabla_x \phi(X,x),x-x^* \>\\
=\<C-\A^*(\lambda),X-X^*\>+\<c-B^\top \lambda,x-x^*\>.
\end{multline}
From some simple algebraic manipulations, we get the following two formulas
\begin{eqnarray}\label{Mar_7_4}
&&\quad \hR\hR^\top-\hR^*\hR^{* \top}=\hR\big(\hR-\hR^*\big)^\top
+\big(\hR-\hR^*\big)\hR^\top-\big(\hR-\widehat{R^*}\big)\big(\hR-\hR^*\big)^\top
\\
\label{Mar_7_5}
&&\quad y\circ y-y^*\circ y^*=y\circ (y-y^*)+( y-y^*)\circ y-(y-y^*)\circ (y-y^*).
\end{eqnarray}
Substituting \eqref{Mar_7_4} and \eqref{Mar_7_5} into \eqref{Mar_7_3} and using the fact that the first $k$ rows of $\hR-\hR^*$ are zero vectors, we get 
\begin{eqnarray}\label{Mar_7_6}
&& \<\nabla _X \phi(X,x),X-X^*\>+\< \nabla_x \phi(X,x),x-x^* \>\\
&=&2\<J_{n-k,n}\(C-\A^*(\lambda)\)\hR,R-R^* \>-\<C_{22}-\tA^*(\lambda),(R-R^*)(R-R^*)^\top\>\nonumber \\
&&+
\,2\,\< \(c-B^\top \lambda\)\circ y,y-y^* \>-\< \(c-B^\top \lambda\), \(y-y^*\)\circ \(y-y^*\)\>
\nonumber \\
&=&\< \g f_r(R,y), (R,y)-(R^*,y^*)\>-\<C_{22}-\tA^*(\lambda),(R-R^*)(R-R^*)^\top\>
\nonumber \\
&& -\< \(c-B^\top \lambda\), \(y-y^*\)\circ \(y-y^*\)\>. \nonumber
\end{eqnarray}
From the Cauchy-Schwarz inequality, we have that
\begin{equation}\label{Mar_7_7}
\< \g f_r(R,y), (R,y)-(R^*,y^*)\>\leq \| \g f_r(R,y) \|\sqrt{\|R-R^*\|^2+\|y-y^*\|^2}.
\end{equation}
Now, apply \eqref{Mar_7_8} and \eqref{Mar_7_7} in \eqref{Mar_7_6}, we get
\begin{multline}\label{Mar_7_9}
\<\nabla _X \phi(X,x),X-X^*\>+\< \nabla_x \phi(X,x),x-x^* \>\\
\leq \| \g f_r(R,y) \|\sqrt{\|R-R^*\|^2+\|y-y^*\|^2}+\sigma_h\|R-R^*\|^2+\sigma_h\|y-y^*\|^2.
\end{multline}
From (ii), we have the following inequality
\begin{multline}\label{Mar_8_1}
\alpha\( \|X-X^*\|^2+\|x-x^*\|^2 \)/2\\
\geq \alpha\sigma_r(\hR^*)^2\|R-R^*\|^2/4+ \alpha \sigma_{\min}(y^*)^2\|y-y^*\|^2/4.
\end{multline}
Combining \eqref{Mar_7_1}, \eqref{Mar_7_9} and \eqref{Mar_8_1}, we have that
\begin{multline}\label{Mar_8_2}
\| \g f_r(R,y) \|\sqrt{\|R-R^*\|^2+\|y-y^*\|^2}\\
\geq \big( \alpha \sigma_r(\hR^*)^2/4-\sigma_h\big)\|R-R^*\|^2+\(\sigma_{\min}(y^*)^2 /4-\sigma_h \)\|y-y^*\|^2.
\end{multline}
From \eqref{Mar_8_2} and the definition of $\sigma_h$ in \eqref{Mar_9_1}, we have that 
\begin{equation}\label{Mar_9_3}
\| \g f_r(R,y) \|\sqrt{\|R-R^*\|^2+\|y-y^*\|^2}\geq \sigma_h \(\|R-R^*\|^2+\|y-y^*\|^2\).
\end{equation}
Note that in deriving the above inequality, when $\|y^*\|=0,$ we have that $\|y\|=0$ because ${\rm supp}(y^*)={\rm supp}(y).$ In this case, the term $\|y-y^*\|^2$ can be ignored. From \eqref{Mar_9_3}, we obtain that
\begin{equation}\label{Mar_8_3}
\| \g f_r(R,y) \|\geq \sigma_h\sqrt{\|R-R^*\|^2+\|y-y^*\|^2},
\end{equation}
which can be derived from dividing both sides of \eqref{Mar_8_2} by $\sqrt{\|R-R^*\|^2+\|y-y^*\|^2}$ when it is nonzero. Note that when $\sqrt{\|R-R^*\|^2+\|y-y^*\|^2}=0,$ then $(R,y)=(R^*,y^*)$ and both sides of \eqref{Mar_8_3} are zero. From \eqref{Mar_9_2} and \eqref{Mar_8_3}, we know that \eqref{PLineq} holds with the parameter $\eta= \sigma_h^2/2\kappa>0$.
\end{proof}

With Lemma~\ref{PLlem}, we can directly prove the following theorem about the local linear convergence of the Riemannian gradient descent method in Algorithm~\ref{alg1}.

\begin{theorem}\label{Thmlinear}
Suppose Assumption~\ref{asslicq} holds, $\bar{r}\leq r\leq n$ and $\{\(R_k,y_k\)\}_{k\in \mathbb{N}}$ is the iteration sequence generated from Algorithm~\ref{alg1} for solving (SDP). Suppose further that $\g f_r$ is locally Lipschitz on $\BA.$ Moreover, assume that the following conditions hold
\begin{itemize}
\item [(i)] $\(R_k,y_k\)$ converges to some $(0,0)-$stationary point $(R^*,y^*)\in \BA$ that satisfies Assumption~\ref{rquad} with parameters $\alpha,\epsilon>0$.
\item [(ii)] There exists $N\in \mathbb{N}_+$ such that for any $k\geq N,$ $\epsilon_i=0,$ $\(R_k,y_k\)\in \BA$ and ${\rm supp}(y_k)={\rm supp}(y).$ Also, the sufficient descent condition \eqref{Sep_25_1} holds.
\end{itemize}
Then $f_r(R_k,y_k)$ converges to $f_r(R^*,y^*)$ R-linearly.
\end{theorem}

\begin{proof}
Because $\(R_k,y_k\)$ converges to $(R^*,y^*)$ and $(R_k,y_k)\in \BA$ when $k\geq N,$ we see that there exists $N_1\geq N$ such that for any $k\geq N_1,$ $\max\{ \|R_k-R^*\|,\|y_k-y^*\| \}<\epsilon$ and $\| \g f(R_k,y_k) \|\leq 1.$ From (ii), \eqref{Sep_25_1} and Lemma~\ref{PLlem}, we have that the following inequality holds
\begin{equation}\label{Mar_8_5}
f_r(R_{k+1},y_{k+1})-f_r(R_k,y_k)\leq -\gamma_g\eta\( f_r(R_k,y_k)-f^* \),
\end{equation}
where $f^*=f_r(R^*,y^*)$ and we have used $\| \g f_r(R_k,y_k) \|\leq 1$ in the right-hand side of \eqref{Sep_25_1}. Rearranging \eqref{Mar_8_5}, we get the following inequality
\begin{equation}\label{Mar_8_6}
f_r(R_{k+1},y_{k+1})-f^*\leq \(1-\gamma_g\eta\)\( f_r(R_k,y_k)-f^* \).
\end{equation}
From \eqref{Mar_8_6}, we  see that if $\gamma_g\eta\geq 1,$ then $f_r(R_{k+1},y_{k+1})=f^*,$ so the algorithm terminates in finite number of 
iterations which is a special case of the R-linear convergence. If $\gamma_g\eta< 1,$ then \eqref{Mar_8_6} is exactly the R-linear convergence inequality with parameter $\(1-\gamma_g\eta\).$
\end{proof}

\begin{remark}\label{remlinear}
Condition (ii) requires us to identify the correct rank and support of the variables $R$ and $y$, respectively. This is indeed necessary because as shown in \cite{zhang2021preconditioned}, the gradient descent only has sublinear convergence for over-parameterized quadratic SDP problems. \end{remark}

\section{Generic smoothness and random perturbation strategy}\label{Sec-randper}

The algorithmic design and theoretical analysis in the previous section are based on the LICQ property, i.e., Assumption~\ref{asslicq}. However, when (SDP) is degenerate, the singularity of the linear system \eqref{projeq} will not only slow down the computation of the projection and retraction but also affect the feasibility of the dual variable $\lambda.$ In this section, we will discuss how to deal with the degeneracy issue.

\subsection{Generic smoothness of (SDPR)}
In this subsection, we will show that the LICQ property of (SDPR) will hold as long as some entries of $b$ are generic. Our main tool is the following Morse-Sard Theorem (see 1.3 of Chapter 3 in \cite{DT}). 

\begin{theorem}[Morse-Sard Theorem]\label{MST}
Let $\M,\N$ be manifolds of dimension $m,n$ and $h:\M\rightarrow \N$ a $C^r$ mapping. If $r>\max\{0,m-n\},$ 
then $h\(\Sigma_h\)$ has measure zero in $\N,$ where $\Sigma_h$ is the set of critical points of $h.$
\end{theorem}

In the above theorem, a point $x\in \M$ is a critical point of $h$ if its differential (also known as pushforward) ${\rm T}_x h:{\rm T}_x\M\rightarrow {\rm T}_{h(x)}\N$ is not surjective. Otherwise, $x$ is called a regular point. Here ${\rm T}_x\M$ and ${\rm T}_{h(x)}\N$ are the tangent spaces of $\M$ and $\N$ at $x$ and $h(x)$ respectively\footnote{For definitions in manifolds and mappings, please refer to chapter 1 of \cite{DT}.}. With Theorem~\ref{MST}, we state our result.

\begin{theorem}\label{randper}
Consider the following set 
$$\N(b,d):=\left\{ x\in \R^n:\ \F(x)=b \right\}\cap \left\{ x\in \R^n:\ \G(x)=d  \right\},$$ 
where $b\in \R^p,d\in \R^q$ and $\F:\R^n\rightarrow \R^p,$ $\G:\R^n\rightarrow \R^q$ are smooth mappings. Suppose every point in $\left\{ x\in \R^n:\ \F(x)=b \right\}$ satisfies the LICQ property. Then for a generic $d\in \R^q,$ every point in $\N(b,d)$ satisfies the LICQ property.
\end{theorem}
\begin{proof}
Define $\M:=\left\{ x\in \R^n:\ \F(x)=b \right\}.$ Consider the mapping $\G_\M:\M\rightarrow \R^q,$ which is the restriction of $\G$ on $\M.$ Because every point in $\M$ satisfies LICQ and $\F$ is a smooth mapping, $\M$ is a smooth submanifold embedded in $\R^n.$ Since $\G$ is a smooth mapping in $\R^n,$ its restriction $\G_\M$ is also smooth. For any $x\in \M,$ the tangent space of $\M$ at $x$ is given by the following formula:
\begin{equation}
{\rm T}_x\M=\left\{ h\in \R^n:\ {\rm J}_\F(x)h=0 \right\},
\end{equation}
where ${\rm J}_\F(x)\in \R^{p\times n}$ is the Jacobian of $\F$ at $x.$ It has full row rank because of the LICQ property of $\M$. The differential of $\G_\M$ at $x\in \M$ is a linear mapping ${\rm T}_x\G_\M: {\rm T}_x\M\rightarrow \R^q$ such that for any $h\in {\rm T}_x\M,$ ${\rm T}_x\G_\M(h)={\rm J}_\G(x)h.$
\begin{claim}\label{claimM}
For any $x\in \M,$ $x$ is a critical point of $\G_\M$ if and only if $x$ doesn't satisfy LICQ in $\N(b,\G(x)).$
\end{claim}
\noindent{\bf Proof of claim.} From the definition of a critical point, we know that $x\in \M$ is a critical point if and only if following inequality holds
\begin{equation}\label{Sep_24_1}
\rr\({\rm J}_\G(x) {\rm J}_\F(x)_\perp^\top {\rm J}_\F(x)_\perp\)\leq q-1,
\end{equation}
where ${\rm J}_\F(x)_\perp\in \R^{(n-p)\times n}$ is a row-orthogonal
 matrix whose row space is the orthogonal complement of the row space of ${\rm J}_\F(x).$ We have that
\begin{multline}\label{Sep_24_2}
\rr\( {\rm J}_\G(x) {\rm J}_\F(x)_\perp^\top {\rm J}_\F(x)_\perp \)+\rr\( {\rm J}_\F(x) \)=\rr\( \begin{bmatrix} {\rm J}_\G(x) {\rm J}_\F(x)_\perp^\top {\rm J}_\F(x)_\perp \\ {\rm J}_\F(x) \end{bmatrix} \)\\
=\rr\( \begin{bmatrix} {\rm J}_\G(x) {\rm J}_\F(x)_\perp^\top {\rm J}_\F(x)_\perp+{\rm J}_\G(x) \(I_n-{\rm J}_\F(x)_\perp^\top {\rm J}_\F(x)_\perp\) \\ {\rm J}_\F(x) \end{bmatrix} \) =\rr\( \begin{bmatrix} {\rm J}_\G(x) \\ {\rm J}_\F(x) \end{bmatrix} \), 
\end{multline}
where the second equality comes from the fact that $I_n-{\rm J}_\F(x)_\perp^\top {\rm J}_\F(x)_\perp$ is the projection matrix of the row space of ${\rm J}_\F(x).$ Because ${\rm J}_\F(x)$ has full row rank, $\rr\( {\rm J}_\F(x) \)=p.$ Therefore, \eqref{Sep_24_1} is equivalent to $\rr\( \begin{bmatrix} {\rm J}_\G(x) \,;\, {\rm J}_\F(x) \end{bmatrix} \)\leq p+q-1,$ which means that $x$ doesn't satisfy LICQ in $\N(b,\G(x)).$ 

From Claim~\ref{claimM}, we know that the set of $d\in \R^q$ such that $\N(b,d)$ contains a singular point is exactly $\G_\M\(\Sigma_{\G_\M}\),$ where $\Sigma_{\G_\M}$ is the set of critical points of $\G_\M.$ From Theorem~\ref{MST}, $\G\(\Sigma_{\G_\M}\)$ has measure zero in $\R^q.$ Thus, for a generic $d\in \R^q,$ except for a zero-measure set, every point in $\N(b,d)$ satisfy LICQ.
\end{proof}


Theorem~\ref{randper} shows that the smoothness of (SDPR) is a generic property. It can also be used to verify that for linear SDP problems, the quadratic growth Assumption~\ref{rquad} is a generic property, as shown in the next proposition. Note that our
result is more general than that established in \cite{genesmo} in that we do not 
need to assume that the linear mapping $\A$ is generic.

\begin{proposition}\label{genegth}
Consider the standard linear SDP problem \eqref{SDP0}. For a generic $b\in \R^m$ and $C\in \S^n,$ the primal and dual constraint nondegeneracy hold everywhere.
As a result, the quadratic growth condition \eqref{rquad} holds at the optimal solution.
\end{proposition}

\begin{proof}
The dual problem and its factorized version of \eqref{SDP0} is as follows:
\begin{eqnarray}\label{SDPD0}
&\max\left\{b^\top y:\ C-\A^*(y)\in \S^n_+,\ y\in \R^m\right\},
\\
\label{SDPD0F}
&\max\left\{b^\top y:\ C-\A^*(y)=YY^\top,\ y\in \R^m,\ Y\in \R^{n\times n}\right\}.
\end{eqnarray}
From Theorem~\ref{randper}, we know that for a generic $b$ and $C$, the LICQ properties hold everywhere in \eqref{SDPLR} and \eqref{SDPD0F} for every $r\geq 1.$ Because it has been shown in \cite[Prop. 6.2]{Boumal3} and \cite[Section 4]{lourencco2018optimality} that the primal and dual constraint nondegeneracy is equivalent to the LICQ properties of \eqref{SDPLR} and \eqref{SDPD0F} respectively, we know that the primal and dual constraint nondegeneracy also hold everywhere in \eqref{SDP0} and \eqref{SDPD0}. Note that the primal nondegeneracy implies the uniqueness of the dual variable, which is exactly the multiplier coming from \eqref{projeq}. From Remark~\ref{remquad}, we know that the quadratic growth condition holds at the optimal solution.
\end{proof}

\subsection{Random perturbation and adaptive preconditioning}\label{percon}
Suppose the constraints in (SDPR) can be separated into two parts $[m]=\I\sqcup \J$ such that the LICQ of $\{(R,y)\in \R^{(n-k)\times r}\times \R^p:\ \K_\I(\hR\hR^\top,y\circ y)=b_\I\}$ is always satisfied, where $b_\I$ denotes the components of $b$ with indices in $\I$; 
a similar definition is used for $\K_\I$.
Then from Theorem~\ref{randper}, after adding a random perturbation with normal distribution to $b_{\J},$ the LICQ will hold everywhere with probability 1. This implies that we don't have to perturb all the entries of $b$ to achieve smoothness. After adding a random perturbation $v\in \R^m$  such that $\|v\|=\epsilon$ to $b$, we are solving a new problem (SDPR$_\epsilon$), i.e., replacing $b$ by $b_\epsilon := b+v$ in (SDPR). Let $(R_\epsilon,y_\epsilon)$ be an $(\epsilon_g,\epsilon_h)-$stationary point of (SDPR$_\epsilon$) and $\lambda_\epsilon$ be the unique solution of $\eqref{projeq}.$ Then $(R_\epsilon,y_\epsilon,\lambda_\epsilon)$ satisfies conditions \eqref{appsta1} and \eqref{appsta2} because they are not related to $b$ explicitly. Moreover, we have that 
$\big\| \K\big(\hR_\epsilon\hR_\epsilon^\top,y_\epsilon\circ y_\epsilon\big)-b
\big\|\leq \big\| \K\big(\hR_\epsilon\hR_\epsilon^\top,y_\epsilon\circ y_\epsilon\big)-(b+v)\big\|+\| v\|=\epsilon.$ This means that the primal feasibility of $(R_\epsilon,y_\epsilon)$ in (SDPR) can be directly controlled by the size of the random perturbation. Therefore, we can solve (SDPR$_\epsilon$) to get an approximate KKT solution of the original problem (SDPR). 
\begin{remark}\label{slater}
Note that after adding a random perturbation to $b,$ the feasible set might become empty. However, we can ensure that the perturbed feasible set is non-empty if the perturbation is sufficiently small, provided (1) the linear mapping $\H: (X,x)\rightarrow \(\K(X,x),X_{1:k,1:k}\)$ is surjective and (2) Slater's condition holds for (\ref{SDP}). The reason is that when the Slater's condition holds, there will exist a feasible solution 
$(\tilde{X},\tilde{x})$ such that $\tilde{X}$ is positive definite and $\tilde{x}>0.$ In this case, there exists $\delta_1>0$ such that such that ${\rm B}_{\delta_1}\(\(X,x\)\)\subset \S^n_+\cap \R^p_+.$ Because $(X,x)\rightarrow (\K(X,x),X_{1:k,1:k})$ is a surjective linear mapping and $\H(\tilde{X},\tilde{x})=\(b,I_k\)$, we have that there exists $\delta_2>0$ such that $\({\rm B}_{\delta_2}(b),I_k\)\subset \H( {\rm B}_{\delta_1}((\tilde{X},\tilde{x})) ).$ This means as long as the size of the perturbation $\epsilon$ is less than $\delta_2,$ the feasible set of SDPR$_\epsilon$ 
is non-empty.
\end{remark}

In practice, the linear system~\eqref{projeq} might still be ill-conditioned after a small random perturbation. To ameliorate the issue, we will
 combine random perturbation with an adaptive preconditioning strategy. The idea is that we use a preconditioned conjugate gradient method (PCG) to solve the linear system \eqref{projeq} and set the maximum iteration of PCG to be some integer $\T_{\rm cg}\geq 2.$ In our implementation, we choose $\T_{\rm cg}$ to be 20 when $m<10000$ and 50 when $m\geq 10000$. When the number of PCG iterations reaches $\T_{\rm cg}$ for the first time, we add a random perturbation to $b$ to ensure smoothness. Before the random perturbation, we use  $\DD(Q_i)^{-1}$ as the preconditioner of PCG. Every time when the number of PCG iterations reaches $\T_{\rm cg},$ we will update the preconditioner $Q_i^{-1}=L^{\top-1}_iL^{-1}_i,$ where $Q_i\in \S^m$ is the coefficient matrix of \eqref{projeq} in the $i$th iteration and $L_iL_i^\top$ is its (sparse) Cholesky factoriztion. By doing so, the number of PCG iterations is upper bounded by $\T_{\rm cg}.$ Moreover, the random perturbation will avoid updating the preconditioner too frequently. This is because from Theorem~\ref{randper}, the feasible set of (SDPR$_\epsilon$) becomes a manifold with probability one. If $(R_i,y_i)\rightarrow (R_\epsilon,y_\epsilon),$ then $Q_i$ will also converge to some positive definite matrix $Q_\epsilon.$ Since $Q_\epsilon$ is positive definite, we know that $Q_i^{-1}\rightarrow Q_\epsilon^{-1}.$ This implies that our preconditioner would not be updated when the iterations is close enough to the optimal solution. After discussing the adaptive preconditioning strategy, we are now able to study the computational complexity of Algorithm~\ref{alg1}.
 
\subsection{Computational complexity of Algorithm~\ref{alg1}}\label{subsec-complx}
In this subsection, we discuss the computational complexity of Algorithm~\ref{alg1}. The main computational cost comes from the projection and retraction mappings, which involve solving only a small number (typically less than $10$) of linear systems like \eqref{projeq} because of the fast (quadratic) convergence of the Newton retraction. We first consider the complexity of solving one linear system without employing a preconditioner. For simplicity of illustration, we assume that there is no vector variable and $k=0.$ In this case, the linear operator on the left-hand side of \eqref{projeq} becomes 
\begin{equation}\label{linop}
2\, \A\( (\A^*(\lambda)R)R^\top +R\(\A^*(\lambda)R\)^\top \).
\end{equation}
Suppose the matrices $A_i$'s have a sparse plus low-rank structure with $A_i=S_i+U_iU_i^\top$ such that $S_i$ is a sparse matrix with $n_i$ nonzero entries and $U_i\in \R^{n\times k_i}$ with small $k_i$. Then computing $\A^*(\lambda)R=\sum_{i=1}^m\lambda_i\(S_i R+U_i\(U_i^\top R\)\)$ has complexity  $\O\(\sum_{i=1}^m n_i r+nrk_i\).$ After we have obtained $V:=\A^*(\lambda)R,$ computing $\A\(VR^\top+RV^\top\)$ also has complexity $\O\( \sum_{i=1}^m n_i r +nrk_i\)\footnote{We only compute the entries of $RR^\top$ and $UV^\top+VU^\top$ when they appear in the sparsity pattern of matrices $S_i.$}.$ When using the conjugate gradient (CG) method to solve \eqref{projeq}, the complexity becomes $\O\( \T\(\sum_{i=1}^m n_i r+nrk_i\) \),$ where $\T$ is the number of CG steps. If we further have  that $\sum_{i=1}^mn_i=\O(m+n),$ $\sum_{i=1}^m k_i=O(1),$ as is often the case for SDP problems coming from real applications, then the computational complexity becomes $\O\( \T(m+n)r \).$ Note that this is just the complexity of solving one linear system. Now we move on study the complexity for Algorithm~\ref{alg1} with the adaptive preconditioning strategy. The linear operator in \eqref{linop} 
is an $m\times m$ symmetric positive definite matrix 
 $Q$ with its $(i,j)$ entry given by $Q_{ij} =4\<A_iR,A_jR\>.$ We have that computing all the entries of $Q$ has the complexity of $\O\( (m+n)^2r+(m+n)r \)\footnote{Note that for each $i\in [m]$ we only have to compute $A_iR$ once.}.$ We should note that the cost for computing $Q$ could be much smaller than $\O\( (m+n)^2r+(m+n)r \)$ when $Q$ itself is sparse, for example, in the case of max-cut problems, the complexity is $\O(nr).$ Let $\T_{\rm alg}$ be the total iterations of Algorithm~\ref{alg1}, $\T_{\rm ch}$ be the number of Cholesky decompositions and $r$ be the upper bound on the rank parameter in Algorithm~\ref{alg1}. Then the total computational complexity becomes 
\begin{equation}\label{totalcomp}
\O\( \T_{\rm alg}\T_{\rm cg}\(m^2+(m+n)r\)+\T_{\rm ch} \(m^3+(m+n)^2r+(m+n)r\) \),  
\end{equation}
where $m^2$ comes from the forward and backward substitutions of the preconditioner $L^{\top-1}L^{-1}$ and $m^3$ comes from the Cholesky decomposition of $Q$ without taking into account of any possible 
sparsity in $Q$. But we should mention that 
for SDP problems with sparse constraint matrices $\{A_i : i\in [m]\},$ it is possible for $Q$ to be sparse, for example,
in the case of Lov\'{a}sz Theta SDPs and max-cut SDPs.
In which case the cost of computing 
the Cholesky decomposition of $Q$ could be 
much smaller than $O(m^3)$ since
we can use sparse Cholesky decomposition.
In contrast, for the interior-point methods in \cite{TTT} and the semi-smooth-Newton augmented Lagrangian
method in \cite{SDPNAL}, the key $m\times m$ matrix arising in each iteration is fully dense even if 
$\{A_i : i\in[m]\}$ are sparse matrices. The possibility of having a sparse $Q$ is one of the key 
computational advantage of our feasible method.

 In order to further make use of the sparsity of $Q$ to get a sparse factor $L$, we further use a hybrid Cholesky decomposition strategy in our
 numerical experiments. In detail, when the number of constraints $m\geq 10000$, we will first try incomplete Cholesky decomposition to see whether it can already control the number of PCG iterations to be less than $\T_{\rm cg}$ and only switch to the exact Cholesky decomposition when the incomplete one does not perform well and exact decomposition becomes necessary.

\section{Dual refinement}\label{Sec-RRdual}
Suppose $(X^*,x^*,\lambda^*)$ is the outcome of Algorithm~\ref{alg1} for solving (SDP) and $(C^*,c^*)=\nabla_{X,x}\phi(X^*,x^*).$ As mentioned in Remark~\ref{remKKT}, the dual variable $\lambda^*$ can be recovered by solving the linear system
in  \eqref{projeq}. The uniqueness of the solution in \eqref{projeq} is equivalent to the LICQ property which can be achieved by our random perturbation strategy. However, the dual feasibility (condition (iii) in Lemma~\ref{Optlem}) of the recovered $\lambda^*$ is sometimes unsatisfactory when the linear system \eqref{projeq} is ill-conditioned. In order to improve the dual feasibility, we consider the following problem:
\begin{multline}\label{dproj}
\min\Big\{  \<C^*,X\>+\<c^*,x\>-\<\lambda^*,\K(X,x)-b\>+\frac{1}{2}\| \K(X,x)-b \|_M^2: 
\\
X_{1:k,1:k}=I_k,\ X\in \S^n_+,\ x\in \R^p_+ \Big\},
\end{multline}
where $M\in \S^n$ is positive definite and $\| . \|_M$ is the weighted norm defined to be $\|A\|_M:=\sqrt{\<A,MA\>}.$ The low rank formulation of problem \eqref{dproj} is as follows:
\begin{eqnarray}\label{lrdproj}
\qquad \min \left\{
\begin{array}{ll}
 \<C^*,\hR\hR^\top\>+\<c^*,y\circ y\>-\<\lambda^*,\K(\hR\hR^\top,y\circ y)-b\> 
 \\[3pt]
 + \frac{1}{2}\| \K(\hR\hR^\top,y\circ y)-b \|_M^2
\end{array}
:
\begin{array}{ll}
  R \in \R^{(n-k)\times r}\\
  y\in \R^p
  \end{array}
 \right\}.
\end{eqnarray}
Since problem \eqref{dproj} is a quadratic SDP, it is equivalent to problem \eqref{lrdproj} as long as $r\geq \lceil \sqrt{k(k+1)+2m} \rceil.$ Moreover, \eqref{lrdproj} is an unconstrained problem, which can be directly solved by Algorithm~\ref{alg1}. Suppose $(R,y)$ is an optimal solution of \eqref{lrdproj} satisfying the optimality conditions in Lemma~\ref{Optlem}. Let $\Lambda$ be the dual variable of $X_{1:k,1:k}=I_k$, as derived in \eqref{Aug_8_18}. We have that 
\begin{eqnarray}\label{Jan_4_2}
C^*-\A^*(\lambda^*)+\A^*\(M\big(\K(\hR\hR^{\top},\;y\circ y)-b\big)\)-I_{n,k}\Lambda I_{n,k}^\top\in \S^n_+,
\\
\label{Jan_4_3}
c^*-B^\top(\lambda^*)+B^\top\(M\big(\K(\hR\hR^{\top},\; y\circ y)-b\big)\)\in \R^p_+.
\end{eqnarray}
By choosing $\lambda:=\lambda^*-M\big(\K(\hR\hR^{\top},y\circ y)-b\big)$, we get the dual variable that satisfies the dual feasibility conditions
$C^*-\A^*(\lambda)-I_{n,k}\Lambda I_{n,k}^\top \in S^n_+$ and $c^*-B^\top \lambda\in \R^p_+.$ Note that we solve problem \eqref{dproj} just to improve the dual feasibility of $\lambda^*.$ The primal optimal solution $(R,y)$ of \eqref{lrdproj} won't be used. In practice, we choose $M$ to be $Q^{-1},$ where $Q$ is the coefficient matrix of \eqref{projeq} mentioned in subsection~\ref{percon}. This usually makes problem \eqref{lrdproj} much easier to solve.
\section{Numerical experiments}\label{Sec-numerexp}
In this section, we test the efficiency of Algorithm~\ref{alg1} in comparison to other solvers. In the gradient descent step of Algorithm~\ref{alg1}, we use Riemannian gradient descent with BB step and non-monotone line search (see \cite{RBB1,non-mono,GEP}). In the escaping saddle point stage, we use Armijo line-search with $c_h'=0.5,$ $t_0=1$ and $\delta=\frac{1}{2}$ as in the proof of Proposition~\ref{prophess}. We choose the tolerance $\epsilon_g=10^{-5}, \epsilon_h=10^{-6}$. For some problems, we also choose a smaller $\epsilon_g$ to get a more accurate solution. We set the tolerance of the Newton retraction to be $10^{-8}$ and the tolerance of PCG to be $10^{-9}.$ Suppose $(R,y,\lambda)$ is the output of Algorithm~\ref{alg1}. Define $X:=\hR\hR^\top$ and $x:=y\circ y.$ We use the following KKT residue to measure the accuracy of the solution for (SDP):
\begin{eqnarray}\label{Pfeas}
& {\rm Rp}:=\frac{\sqrt{\|\K\(X,x\)-b\|^2+\|X_{1:k,1:k}-I_k\|^2}}{1+\sqrt{\|b\|^2+k}},\notag
\\
\label{Dfeas}
& {\rm Rd}:=\frac{ \sqrt{\|\Pi_{\S^{n}_-}\({\bf S}\)\|^2+\|\Pi_{\R^p_-}\({\bf s}\)\|^2  }}{1+\sqrt{\|C\|^2+\|c\|^2}},\ {\rm Rc}:=\frac{ |\<{\bf S},X\>+\<{\bf s},x\>| }{1+\sqrt{\|C\|^2+\|c\|^2}}, \notag
\end{eqnarray}
where ${\bf S}\in \S^n$ is defined as in \eqref{Aug_8_17} and ${\bf s}:=c-B^\top\lambda.$ We don't check the constraints $X\in \S^n_+$ and $x\in \R^p_+$ because they are strictly satisfied from their factorization forms. We call the practical implementation of the algorithm SDPF, where ``F'' indicates that it is a feasible method.

For linear SDP problems, we set $r_0$ to be the theoretical rank upper bound $\lceil \sqrt{2m+k(k+1)}\rceil.$ 
We compare SDPF with the solver SDPLR \cite{BM1}\footnote{Source codes from \href{https://sburer.github.io/projects.html}{https://sburer.github.io/projects.html}.} and 
SDPNALplus \cite{SDPNALp2}\footnote{Source codes from \href{https://blog.nus.edu.sg/mattohkc/softwares/sdpnalplus/}{https://blog.nus.edu.sg/mattohkc/softwares/sdpnalplus/}.}. 
The former is a low rank SDP solver based on ALM. The latter
 is a convex SDP solver based on ALM and ADMM. We have also tested SketchyCGAL by Yurtsever et al \cite{CGAL}. However, we find that this solver cannot find an accurate solution and is always quite slow, so we don't show the results of this solver. For quadratic SDP problems, we compare SDPF with QSDPNAL\footnote{Source codes from \href{https://blog.nus.edu.sg/mattohkc/softwares/qsdpnal/}{https://blog.nus.edu.sg/mattohkc/softwares/qsdpnal/}.} \cite{QSDPNAL}. The above algorithms have different computational bottlenecks. For SDPF, the most expensive part comes from solving the linear systems in projection and retraction mappings. For SDPLR, it is costly to solve the ALM subproblem in each iteration especially for degenerate problems. For SDPNAL and QSDPNAL, the computational challenge lies in the full eigenvalue decomposition of dense $n\times n$ matrices
 and solving the fully dense $m\times m$ linear systems in the semi-smooth Newton method 
 for each subproblem.

 In the numerical experiments,  we set the accuracy tolerance to be $10^{-6}$ and maximum running time to be 3600s for all solvers. We don't show the results in the following tables if an algorithm has reached the maximum running time but the outcome is still inaccurate. For linear SDP problems\footnote{The quadratic SDP problem we test is unconstrained for Algorithm~\ref{alg1}.}, apart from the running time, we also display $\T_{\rm alg}$: total iterations, $\T_{\rm cg}$: total CG iterations, $\T_{\rm lin}$: total linear systems solved and $\T_{\rm ch}$: total Cholesky decompositions. If the problems is smooth and preconditioner is not used, we will omit $\T_{\rm ch}$ in the table. For SDPLR, $\T_{\rm alg}$ means the number of outer iterations i.e., number of ALM subproblems solved, $\T_{\rm cg}$ means the total limited memory BFGS steps in solving the ALM subproblems. For SDPNAL, $\T_{\rm alg}$ means the number of ADMM iterations plus the number of semi-smooth Newton subiterations. It is also equal to the number of eigenvalue decompositions taken in this algorithm. For SDPNAL, $\T_{\rm cg}$ means the total number of Krylov solver iterations taken in solving the linear systems in the semi-smooth Newton method, $\T_{\rm lin}$ denotes the number of linear systems that have been solved in the semi-smooth Newton method. All the experiments are run using Matlab R2021b on a workstation with Intel Xeon E5-2680 v3 @ 2.50GHz processor and 128GB RAM.

\subsection{Random sparse SDP problem}

We first test on some randomly generated linear SDP problems following the generating procedure in subsection 5.2 of \cite{RegSDP}. We set the number of nonzero entries of the upper triangular part of   matrix $A_i$ and each row of $B$ to be 3. Note that our results for (SDP) can easily be extended to multi-block SDP problems. We test on
problems with two matrix variables of the same size $n$ and one vector variable of length $p.$ We choose $k=0$ in this experiment. 

\begin{center}
\begin{tiny}
\begin{longtable}{|c|c|cc|c|c|}
\caption{Comparison of  SDPF , SDPNAL and SDPLR for random sparse SDP problems with small $n$ and large $m=\lceil n^2/4\rceil.$ For the meaning of the notation in the last column, see the paragraph just before section 6.1.}
\\
\hline
problem & algorithm & Rp/Rd/Rc& obj & time & $\T_{\rm alg}$/$\T_{\rm cg}$/$\T_{\rm lin}$/$\T_{\rm ch}$ \\ \hline
\endhead
$n$ = 200 & SDPF & 1.6e-09/3.2e-07/3.2e-06 & 3.4087076e+01& 6.49e+00&36/1757/67/ \\
$p$ = 5000 & SDPNAL & 2.6e-07/1.2e-08/2.9e-07 & 3.4087048e+01& 5.71e+00&249/381/ 2/ \\
$m$ = 10000 & SDPLR & 3.7e-08/2.0e-06/1.8e-06 & 3.4087082e+01& 5.29e+00&19/3365/-/ \\
\hline
$n$ = 300 & SDPF & 1.8e-10/2.3e-07/3.8e-06 & 8.7768271e+01& 1.63e+01&46/2303/90/ \\
$p$ = 11250 & SDPNAL & 9.4e-07/5.8e-08/6.0e-08 & 8.7768405e+01& 1.07e+01&265/360/ 2/ \\
$m$ = 22500 & SDPLR & 2.8e-08/4.5e-07/2.9e-06 & 8.7768295e+01& 1.50e+01&19/3933/-/ \\
\hline
$n$ = 500 & SDPF & 3.2e-10/1.1e-07/1.4e-07 & -2.8714214e+02& 6.91e+01&52/2944/107/ \\
$p$ = 31250 & SDPNAL & 9.8e-07/2.5e-08/5.2e-08 & -2.8714233e+02& 2.84e+01&295/380/ 2/ \\
$m$ = 62500 & SDPLR & 1.6e-08/5.8e-07/1.2e-06 & -2.8714211e+02& 9.38e+01&20/4985/-/ \\
\hline
$n$ = 1000 & SDPF & 8.4e-10/3.4e-08/1.6e-06 & 1.9931863e+01& 6.41e+02&90/5185/199/ \\
$p$ = 125000 & SDPNAL & 8.2e-07/9.2e-09/9.6e-09 & 1.9931857e+01& 1.48e+02&303/420/ 2/ \\
$m$ = 250000 & SDPLR & 8.8e-09/2.2e-07/1.7e-07 & 1.9931885e+01& 1.02e+03&20/9018/-/ \\
\hline
\end{longtable}
\end{tiny}
\end{center}
\begin{center}
\begin{tiny}
\begin{longtable}{|c|c|cc|c|c|}
\caption{Comparison of  SDPF , SDPNAL and SDPLR for random sparse SDP problems with large $n$ and
moderate $m=10n.$}
\\
\hline
problem & algorithm & Rp/Rd/Rc& obj & time& $\T_{\rm alg}$/$\T_{\rm cg}$/$\T_{\rm lin}$/$\T_{\rm ch}$ \\ \hline
\endhead
$n$ = 2000 & SDPF & 6.7e-11/5.5e-07/3.5e-06 & 2.3266098e+01& 8.83e+01&54/904/101/ \\
$p$ = 20000 & SDPNAL & 3.9e-08/8.4e-09/5.7e-07 & 2.3265957e+01& 6.42e+02&547/933/ 6/ \\
$m$ = 20000 & SDPLR & 2.7e-08/4.5e-07/7.6e-06 & 2.3266167e+01& 2.95e+02&19/3278/-/ \\
\hline
$n$ = 3000 & SDPF & 5.2e-11/3.7e-07/2.9e-07 & -1.5008769e+01& 1.81e+02&60/997/113/ \\
$p$ = 30000 & SDPNAL & 3.5e-08/5.9e-08/4.5e-07 & -1.5008945e+01& 1.81e+03&468/896/ 8/ \\
$m$ = 30000 & SDPLR & 2.3e-08/4.2e-07/2.2e-06 & -1.5008746e+01& 6.58e+02&20/3558/-/ \\
\hline
$n$ = 5000 & SDPF & 7.9e-10/5.2e-08/4.2e-07 & 1.2900187e+01& 5.59e+02&82/1391/174/ \\
$p$ = 50000 & SDPNAL & 1.3e-02/4.0e-04/3.1e-02 & 2.8875003e+01& 3.62e+03&220/232/ 7/ \\
$m$ = 50000 & SDPLR & 1.7e-08/4.0e-07/2.0e-05 & 1.2900948e+01& 2.48e+03&20/4559/-/ \\
\hline
$n$ = 10000 & SDPF & 2.1e-11/2.2e-07/6.4e-07 & -5.4679253e+01& 1.80e+03&81/1300/162/ \\
$p$ = 100000 & SDPNAL & - & - & - & - \\
$m$ = 100000 & SDPLR & - & - & - & - \\
\hline
 \end{longtable}
\end{tiny}
\end{center}

Table 1 shows the results for problems with small $n$ and large number of constraints. All the three algorithms can solve all the problems to the required accuracy. SDPNAL is faster than the other two solvers. In particular, it is faster than SDPF by factor of $4.4$ times for the last instance.
This is because for problems with small $n$ and large $m,$ full eigenvalue decompositions are affordable and the optimal solutions
of the tested instances have high rank. These features are conducive for SDPNAL to perform very efficiently.
Table 2 contains problems with large $n$ and moderate number of constraints $m=O(n).$ Such problems have low rank property, so SDPF and SDPLR perform better than SDPNAL. Our algorithm is the only one that can solve all the instances to the accuracy of 
 less than $5\times 10^{-6}$. 
 It is also faster than SDPLR. This demonstrates the efficiency of our feasible method as compared to ALM. Note that for SDPF, $\T_{\rm lin}$ is close to or sometimes even smaller than $2\T_{\rm alg}.$ 
This is because since SDPF is a Riemannian optimization method, a new iteration point on the tangent space is usually so close to the manifold that even just one step of Gauss-Newton iteration (in the Newton retraction) can reach the retraction accuracy tolerance of $10^{-8}.$


\subsection{Sensor network localization}

Consider the following sensor network localization with square loss for $p$ sensors and $q$ anchors:
\begin{equation}\label{SNL}
\min_{u_1,\ldots,u_\ell\in \R^d}\Big\{ \frac{1}{2}\sum_{(i,j)\in \N}\( \|u_i-u_j\|^2-d_{ij}^2 \)^2+\frac{1}{2}\sum_{(i,k)\in \M}\( \|u_i-a_k\|^2-d_{ik}^2 \)^2  \Big\},
\end{equation}
For $i,j\in [p],$ define $e_{ij}:=e_i-e_j\in \R^n.$ Let $g_{ik}:=[-a_k;e_i]$ for $(i,k)\in \M$ and let $g_{ij}=[0_{d};e_{ij}]$ for $(i,j)\in \N.$ The SDP relaxation of \eqref{SNL} is the following quadratic SDP problem:
\begin{equation}\label{SNLSDP}
\min_{X\in \S^{n+d}_+,\atop X_{1:d,1:d}=I_d }\left\{ \frac{1}{2}\sum_{(i,j)\in \N}\( g_{ij}^\top Xg_{ij}-d_{ij}^2 \)^2+\frac{1}{2}\sum_{(i,j)\in \M}\( g_{ik}^\top Xg_{ik}-d_{ik}^2 \)^2-\frac{1}{n}\<I_{n+d}-aa^\top,X\>\right\},\notag 
\end{equation}
where $a=[\hat{a};\hat{e}]$ with $\hat{e}=e/\sqrt{p+q}$ and $\hat{a}=\sum_{k=1}^q a_k/\sqrt{p+q}.$
The last term in the objective is a  regularization used in \cite{SNL} to prevent the predicted data points from being too close to each other. We choose $d=3,$ $q=8$ and the locations of the eight anchors $\{a_i\}_{i\in [8]}$ are $[\pm 0.3;\pm 0.3;\pm 0.3].$ For the distance $d_{ij}$ and $d_{ik},$ we first randomly generate $p$ data points in $[-0.5,0.5]^3$ using \texttt{rand(n,3)-0.5}. Let $x_i$ be the position of the $i$th sensor. For any $i,j\in [p]$ such that $i\neq j,$ let $(i,j)\in \N$ if $\|x_i-x_j\|\leq \(15/(\pi p)\)^{1/3}.$ 
We choose this ratio so that the average degree of $\N$ is close to $20.$ Let $d_{ij}=(1+0.1\sigma_{ij})\|x_i-x_j\|.$ For any $(i,j)\in [p]\times [q],$ let $(i,j)\in \M$ if $\|x_i-a_j\|\leq 0.3.$ Let $d_{ik}=(1+0.1\sigma_{ik})\|x_i-a_k\|.$ Here $\sigma_{ij}$ and $\sigma_{ik}$ are noise satisfying the standard normal distribution. 
For Algorithm~\ref{alg1}, we set $\epsilon_g=\epsilon_h=10^{-6},$ $r_0=3$ and $\tau=20.$

\begin{figure}
\centerline{\includegraphics[height=2cm,width=8cm]{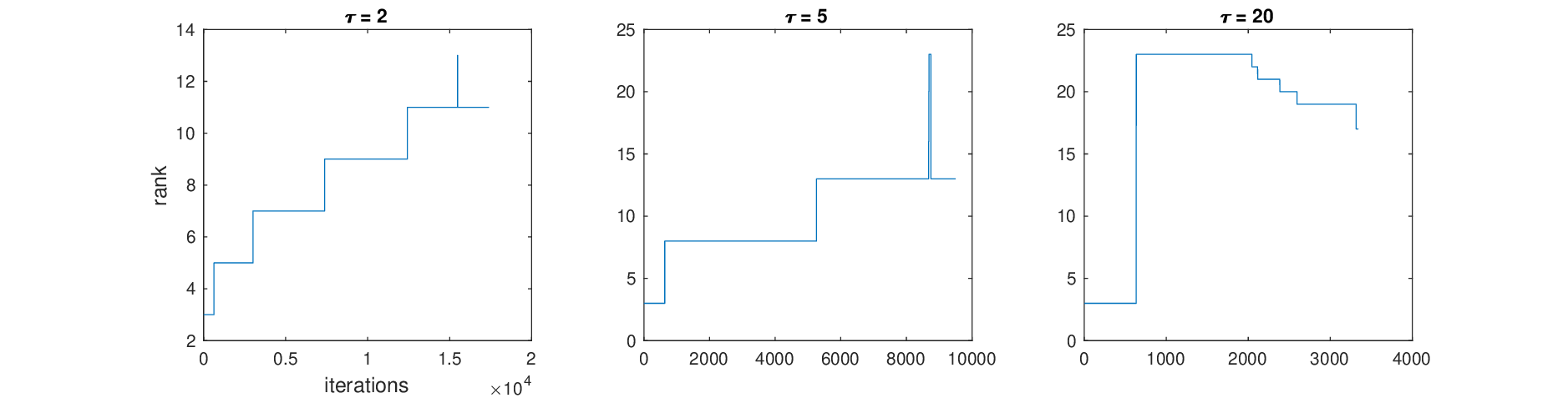}}
\caption{\tiny Search rank of SDPF for the instance $n=1000$ in Table 3 with different rank increasing parameter $\tau.$} 
\label{figrank}
\end{figure}

\begin{center}
\begin{tiny}
\begin{longtable}{|c|c|cccc|l|}
\caption{Comparison of SDPF and QSDPNAL for sensor network localization SDP. The rank 
reported is the numerical rank of the output $R$ of SDPF. }
\\
\hline
problem & algorithm & Rp & Rd & Rc& obj & time \\ \hline
\endhead
n = 300 & SDPF & 0.00e+00 & 4.64e-07 & 1.76e-09 & -2.1939425e-01& 7.03e+00 \\
rank = 12 & QSDPNAL & 4.05e-11 & 6.39e-06 & 2.63e-04 & -2.1937102e-01& 1.38e+02 \\
\hline
n = 400 & SDPF & 0.00e+00 & 3.04e-07 & 6.66e-09 & -2.1655246e-01& 5.58e+00 \\
rank = 13 & QSDPNAL & 3.93e-11 & 3.38e-07 & 3.55e-04 & -2.1651659e-01& 2.29e+02 \\
\hline
n = 500 & SDPF & 0.00e+00 & 1.24e-07 & 2.82e-08 & -2.1084873e-01& 9.70e+00 \\
rank = 15 & QSDPNAL & 4.43e-11 & 2.24e-06 & 1.29e-04 & -2.1075949e-01& 4.62e+02 \\
\hline
n = 1000 & SDPF & 0.00e+00 & 3.12e-07 & 7.38e-09 & -2.1248491e-01& 2.83e+01 \\
rank = 17 & QSDPNAL & - & - & - & -& - \\
\hline
n = 2000 & SDPF & 0.00e+00 & 5.32e-07 & 4.93e-08 & -2.0585246e-01& 7.25e+01 \\
rank = 23 & QSDPNAL & - & - & - & -& - \\
\hline
n = 3000 & SDPF & 0.00e+00 & 8.87e-07 & 2.12e-08 & -2.0396304e-01& 1.21e+02 \\
rank = 23 & QSDPNAL & - & - & - & -& - \\
\hline
n = 5000 & SDPF & 0.00e+00 & 7.45e-07 & 8.56e-09 & -1.8536891e-01& 2.93e+02 \\
rank = 23 & QSDPNAL & - & - & - & -& - \\
\hline
n = 10000 & SDPF & 0.00e+00 & 6.51e-07 & 3.07e-08 & -1.5063329e-01& 6.79e+02 \\
rank = 23 & QSDPNAL & - & - & - & -& - \\
\hline
\end{longtable}
\end{tiny}
\end{center}

From Table 3, we can see that SDPF can solve all the problems to the required accuracy and it is much more efficient than QSDPNAL. The rank of every optimal solution is greater than $r_0=3,$ but our rank-adaptive method can progressively increase the rank and converge to the global optimal solution. Figure~\ref{figrank} shows the rank updating process for different increasing rank parameter $\tau$.


\subsection{Box-constrained quadratic programming}
Vandenbussche and Nemhauser considered the following box-constrained quadratic programming problem in \cite{BQP2}:
\begin{equation}\label{BQP}
\min\left\{ \frac{1}{2}z^\top Q z+a^\top z:\ 0\leq z\leq  e,\ z\in \R^n \right\},
\end{equation}
where $Q\in \S^n,$ $a\in \R^n$ are given data. We consider its SDP relaxation below:
\begin{equation}\label{BQPSDP}
\min\left\{ \frac{1}{2}\<Q,Z\>+a^\top z:\ \dd(Z)+x=z,\ X:=\begin{pmatrix} 1& z^\top \\ z& Z \end{pmatrix}\in \S^{n+1}_+,\ x\in \R^{n}_+ \right\}.
\end{equation}
Observe that \eqref{BQPSDP} has the form of (SDP) with $k=1.$ We generate the data $Q$ and $a$ following the procedure mentioned in Section 4 of \cite{BQP2}.

\begin{center}
\begin{tiny}
\begin{longtable}{|c|c|cc|c|c|}
\caption{Comparison of SDPF, SDPNAL and SDPLR for 
SDP relaxation of box-constrained QP. In the table, $\rho$ is the density of $Q.$}
\\
\hline
problem & algorithm & Rp/Rd/Rc& obj & time & $\T_{\rm alg}$/$\T_{\rm cg}$/$\T_{\rm lin}$/$\T_{\rm ch}$\\ \hline
\endhead
n = 1000 & SDPF & 3.2e-09/1.5e-07/9.0e-09 & -1.4909488e+05& 1.59e+00&186/362/362/ \\
$\rho$ = 2.50e-01 & SDPNAL & 2.5e-07/2.3e-08/1.3e-04 & -1.4909487e+05& 1.10e+03&7301/2769/501/ \\
& SDPLR & 8.3e-07/2.8e-06/7.3e-06 & -1.4909507e+05& 3.08e+01&18/8979/-/ \\
\hline
n = 1000 & SDPF & 4.3e-12/4.3e-08/2.1e-08 & -2.5806628e+05& 1.50e+00&184/353/353/ \\
$\rho$ = 7.50e-01 & SDPNAL & 4.1e-07/4.1e-07/6.1e-05 &- 2.5806644e+05&9.53e+02&5688/2848/488/ \\
& SDPLR & 9.4e-07/2.2e-06/7.3e-06 & -2.5806657e+05& 1.87e+01&19/5283/-/ \\
\hline
n = 2000 & SDPF & 5.8e-09/4.1e-08/1.4e-07 & -4.3135148e+05& 3.81e+00&185/372/372/ \\
$\rho$ = 2.50e-01 & SDPNAL & - & -& -&- \\
& SDPLR & 9.9e-07/7.1e-07/3.4e-05 & -4.3135227e+05& 4.76e+02&18/21690/-/ \\
\hline
n = 2000 & SDPF & 2.7e-11/8.2e-08/9.3e-08 & -7.5485516e+05& 4.19e+00&261/506/506/ \\
$\rho$ = 7.50e-01 & SDPNAL & - & -& -&- \\
& SDPLR & 9.6e-07/6.0e-07/1.0e-04 & -7.5485379e+05& 3.86e+02&20/17561/-/ \\
\hline
n = 5000 & SDPF & 1.2e-09/4.6e-08/1.2e-08 & -1.7227960e+06& 1.96e+01&352/663/663/ \\
$\rho$ = 2.50e-01 & SDPNAL & - & -& -&- \\
& SDPLR & 2.9e-05/1.0e-06/6.3e-04 & -1.7228807e+06& 3.61e+03&18/25940/-/ \\
\hline
n = 5000 & SDPF & 2.3e-10/1.2e-08/3.7e-07 & -3.0025339e+06& 1.97e+01&317/624/624/ \\
$\rho$ = 7.50e-01 & SDPNAL & - & -& -&- \\
& SDPLR & 5.2e-06/3.7e-06/2.2e-03 & -3.0025039e+06& 3.61e+03&18/25380/-/ \\
\hline
n = 10000 & SDPF & 6.6e-10/2.0e-08/2.8e-07 & -4.8841142e+06& 9.46e+01&567/1036/1036/ \\
$\rho$ = 2.50e-01 & SDPNAL & - & -& -&- \\
 & SDPLR & - & - & - & - \\
\hline
n = 10000 & SDPF & 6.7e-09/1.4e-08/1.9e-08 & -8.4655873e+06& 1.08e+02&729/1266/1266/ \\
$\rho$ = 7.50e-01 & SDPNAL & - & -& -&- \\
& SDPLR & - & - & - & - \\
\hline
\end{longtable}
\end{tiny}
\end{center}

From Table 4, we can see that SDPF is significantly faster than SDPLR and SDPNAL. For some instances, SDPF can be more than 100 times faster than SDPLR and nearly 1000 times faster
than SDPNAL. Moreover, SDPF can solve all the problems more accurately than SDPLR and
SDPNAL.
The feasible set of (SDPR) is the Cartesian product of unit spheres, i.e., an oblique manifold. In this case, the coefficient matrix $Q_i$ of the linear system 
in the projection and retraction mapping is diagonal. As mentioned in subsection~\ref{randper}, we use $\DD(Q_i)^{-1}=Q_i^{-1}$ to be the preconditioner before the random perturbation. Thus, our PCG converges in one iteration, which explains why $\T_{\rm cg}=\T_{\rm lin}$ for SDPF. Without dealing with the constraint $X_{11}=1$ as in (SDPR) explicitly, the Jacobian matrix of the constraints can be ill-conditioned, which results in the large number of iterations $\T_{\rm cg}$ for SDPLR to solve the ALM subproblems. This 
demonstrates the effectiveness of our model (SDPR) as compared to the traditional BM factorization.

\subsection{The Lov\'asz Theta SDP}
The problems we have tested before are mostly nondegenerate. In order to verify the effectiveness of our random perturbation strategy, we test the following Lov\'asz Theta SDP problem \cite{Theta} and frequently used as benchmarks for SDP solvers:
\begin{equation}\label{Thetap}
\min\left\{ -\<ee^\top,X\>:\ \<I,X\>=1,\ \A_G\(X\)=0,\ X\in \S^n_+ \right\},
\end{equation}
where $G=(V,E)$ is a simple undirected graph, $\A_G:\S^n\rightarrow \R^m$ is defined by
 $\A_G(X):=\(X_{ij}\)_{ij\in E}.$ We consider two types of graphs in our experiments. The first type of graph, which contains fewer than 3000 vertices, comes from coding theory \cite{Coding}. Such a dataset has already been tested in \cite{SDPNAL}, where SDPNAL achieves good performance. The second group of graphs are large sparse graphs coming from Gset\footnote{Dataset from \href{https://web.stanford.edu/~yyye/yyye/Gset/}{https://web.stanford.edu/~yyye/yyye/Gset/}.} such that $n\geq 5000.$ Note that some of the Gset graphs such as g58 and g59 have the same structure with the only difference in the edge weight. We only show the result of one of those graphs. All the problems we test are degenerate, so we apply our random perturbation strategy mentioned in subsection~\ref{percon}. The KKT residue reported in the following tables is for the original unperturbed problem. 
Our random perturbation is kept to be sufficiently 
small so that the primal residue of the output of SDPF can reach 
the required accuracy of $10^{-6}.$ For some extremely ill-conditioned problems such as g58 and g63, we add a larger perturbation and set the tolerance to be $10^{-5}.$ Apart from the random perturbation, we also apply the dual refinement in Section~\ref{Sec-RRdual} to improve the dual KKT residue of the output. The time spent in the dual refinement is part of the running time.

\begin{center}
\begin{tiny}
\begin{longtable}{|c|c|cc|c|c|}
\caption{Comparison of SDPF, SDPNAL and SDPLR for Lov\'asz Theta SDP in coding theory. In the first column, $n=|V|$ and $m=|E|.$ }
\\
\hline
problem & algorithm & Rp/Rd/Rc& obj & time & $\T_{\rm alg}$/$\T_{\rm cg}$/$\T_{\rm lin}$/$\T_{\rm ch}$  \\ \hline
\endhead
1dc.1024 & SDPF & 7.5e-07/3.8e-07/7.9e-07 & -9.5991877e+01& 5.85e+02&598/25788/1012/9 \\
n = 1024 & SDPNAL & 6.5e-07/6.6e-07/2.9e-08 & -9.5986248e+01& 1.02e+02&251/5165/51/- \\
m = 24063 & SDPLR & 1.0e-06/1.3e-06/2.6e-08 & -9.5987268e+01& 8.68e+02&19/87756/-/- \\
\hline
1dc.2048 & SDPF & 7.5e-07/2.3e-07/4.8e-06 & -1.7474015e+02& 1.68e+03&428/19948/671/42  \\
n = 2048 & SDPNAL & 3.7e-07/3.0e-07/5.8e-09 & -1.7473095e+02& 8.07e+02&273/8668/73/- \\
m = 58367 & SDPLR &  - & - & -& - \\
\hline
1et.1024 & SDPF & 7.6e-07/6.6e-07/8.1e-07 & -1.8423771e+02& 4.92e+01&176/3045/265/9 \\
n = 1024 & SDPNAL & 4.3e-07/4.6e-07/4.9e-09 & -1.8422716e+02& 7.81e+02&4167/18570/167/- \\
m = 9600 & SDPLR & 1.0e-06/5.5e-06/9.6e-09 & -1.8423119e+02& 1.85e+03&19/257078/-/- \\
 \hline
1et.2048 & SDPF & 7.6e-07/8.5e-07/2.3e-07 & -3.4206073e+02& 1.50e+02&218/4998/324/4 \\
n = 2048 & SDPNAL & 2.0e-07/3.2e-07/3.6e-08 & -3.4202941e+02& 2.39e+03&714/26604/214/- \\
m = 22528 & SDPLR & 9.7e-04/8.2e-04/2.9e-05 & -3.5833285e+02& 3.60e+03&14/87907/-/- \\
\hline
1tc.1024 & SDPF & 7.5e-07/5.4e-07/2.6e-07 & -2.0631684e+02& 6.32e+01& 190/2914/312/16 \\
n = 1024 & SDPNAL & 7.0e-07/3.8e-07/1.4e-08 & -2.0630546e+02& 5.45e+02&1615/26745/229/- \\
m = 7936 & SDPLR & 1.0e-06/5.2e-06/4.1e-07 & -2.0630913e+02& 1.18e+03&18/182339/-/- \\
\hline
1tc.2048 & SDPF & 7.6e-07/5.1e-07/7.0e-07 & -3.7467663e+02& 1.78e+02&224/6520/324/8 \\
n = 2048 & SDPNAL & 8.5e-07/2.6e-07/5.6e-07 & -3.7464643e+02& 2.75e+03&1188/30071/269/- \\
m = 18944 & SDPLR & 1.5e-05/4.1e-04/2.4e-04 & -3.7486924e+02& 3.60e+03&15/91296/-/- \\
\hline
1zc.1024 & SDPF & 7.5e-07/8.5e-07/4.6e-07 & -1.2867552e+02& 2.96e+02&430/18411/592/106 \\
n = 1024 & SDPNAL & 6.3e-08/1.5e-07/7.6e-11 & -1.2866668e+02& 4.26e+01&217/710/17/- \\
m = 16640 & SDPLR & 1.0e-06/8.8e-08/3.0e-09 & -1.2866796e+02& 2.70e+03&18/322387/-/-  \\
\hline
1zc.2048 & SDPF & 7.6e-07/5.0e-07/1.6e-08 & -2.3742455e+02& 1.24e+02&85/2661/124/3 \\
n = 2048 & SDPNAL & 5.4e-08/3.8e-09/7.8e-11 & -2.3739996e+02& 4.01e+02&233/1607/33/- \\
m = 39424 & SDPLR & 9.3e-04/1.8e-05/5.2e-07 & -2.4748487e+02& 3.60e+03&16/69904/-/- \\
\hline
\end{longtable}
\end{tiny}
\end{center}

\begin{center}
\begin{tiny}
\begin{longtable}{|c|cc|c|c|}
\caption{Test SDPF for Lov\'asz Theta SDP in Gset. $n,m$ denotes the number of vertices and edges respectively.}
\\
\hline
problem: $n,m$  & Rp/Rd/Rc& obj & time & $\T_{\rm alg}$/$\T_{\rm cg}$/$\T_{\rm lin}$/$\T_{\rm ch}$ \\ \hline
\endhead
g55: 5000,12498  & 7.5e-07/8.1e-07/1.5e-07  &  -2.3243020e+03& 7.96e+01& 231/3824/262/3 \\
\hline
g57: 5000,10000  & 7.0e-09/1.7e-07/8.5e-11 &  -2.4999998e+03& 3.96e+01 & 355/4162/582/36 \\
\hline
g58: 5000,29570  & 7.5e-06/9.4e-06/7.0e-07 &  -1.7833429e+03& 6.68e+02 & 219/10369/326/93 \\
\hline
g60: 7000,17148  & 7.5e-07/6.7e-07/2.7e-07 & -3.2652524e+03& 1.73e+02 & 245/4809/359/3 \\
\hline
g62: 7000,14000  & 2.2e-08/1.5e-07/2.6e-10 & -3.4999998e+03& 5.61e+01 & 259/4138/467/16\\
\hline
g63: 7000,41459  & 7.5e-06/6.9e-06/3.3e-07 & -2.4946356e+03& 1.02e+03 & 160/6174/238/32\\
\hline
g65: 8000,16000  & 7.5e-10/0.0e+00/2.8e-07 & -3.9999992e+03& 9.98e+01& 490/6450/800/65\\
\hline
g66: 9000,18000  & 1.4e-08/2.0e-07/9.6e-12 & -4.4999995e+03& 8.23e+01 &404/4630/641/56\\
\hline
g67: 10000,20000  & 6.0e-09/0.0e+00/1.2e-07 & -4.9999992e+03& 1.19e+02 &342/4156/539/25\\
\hline
g70: 10000,9999$\ $   & 7.6e-07/8.8e-07/3.0e-08 & -6.0777626e+03& 7.12e+02&158/3686/224/13\\
\hline
g77: 14000,28000  & 7.1e-09/5.8e-07/1.9e-10 & -6.9999967e+03& 2.06e+02& 542/7009/932/60\\
\hline
g81: 20000,40000  & 3.4e-09/6.0e-14/3.5e-08 & -1.0000000e+04& 6.37e+02& 1038/10142/1594/94\\
\hline
\end{longtable}
\end{tiny}
\end{center}

From Table 5, we can see that SDPF can solve most of the instances to the required accuracy $10^{-6}.$ For those instances that SDPF didn't reach the tolerance, the KKT residues are already close to $10^{-6}.$ This verifies the robustness of SDPF. As for the speed, both SDPF and SDPNAL are faster than SDPLR but there is no clear winner between SDPF and SDPNAL. Actually, SDPF can be either ten times faster than SDPNAL (1et.2048, 1tc.2048) or nearly 7 times slower than SDPNAL (1zc.1024). This is in line with our expectation because $n$ is moderate and the degeneracy is not significant for this dataset. However, when we test the large sparse Gset graphs that are highly degenerate, we find that neither SDPNAL nor SDPLR can solve any of the instances up to even low accuracy. Thus, we only show the results for SDPF. From Table 6, we can see that SDPF can solve all the problems to the required accuracy. Note that for g58 and g63, we have increased the tolerance to be $10^{-5}$ due to the extreme ill-conditioning. In Algorithm~\ref{alg1}, the number of Cholesky decompositions is usually at least ten times smaller than the number of linear systems solved 
in Algorithm~\ref{alg1}. Besides, the PCG iterations for Algorithm~\ref{alg1} is much smaller than the number of subiterations for SDPLR. This demonstrates the effectiveness of our random perturbation and dual refinement strategies. For the largest graph g81 with 20000 vertices, SDPF can find the solution within 11 minutes. This demonstrates the efficiency of our algorithm for solving large scale low-rank SDP problems.

\section{Conclusion}\label{Sec-conc}
In this paper, we propose a rank-support adaptive feasible 
method to solve general convex  SDP problems with guaranteed convergence. We apply random perturbation and dual refinement strategies to overcome the numerical difficulty caused by the issue of degeneracy. Our work have significantly extended the application range of feasible method for
solving convex SDP problems.


\bibliographystyle{abbrv}
\bibliography{SDPR}

\begin{thebibliography}{10}

\bibitem{manibook}
P.-A. Absil, R.~Mahony, and R.~Sepulchre.
\newblock Optimization algorithms on matrix manifolds.
\newblock In {\em Optimization Algorithms on Matrix Manifolds}. Princeton
  University Press, 2009.

\bibitem{Pretrac}
P.-A. Absil and J.~Malick.
\newblock Projection-like retractions on matrix manifolds.
\newblock {\em SIAM J. Optimization}, 22(1):135--158, 2012.

\bibitem{genesmo}
F.~Alizadeh, J.-P.~A. Haeberly, and M.~L. Overton.
\newblock Complementarity and nondegeneracy in semidefinite programming.
\newblock {\em Mathematical Programming}, 77(1):111--128, 1997.

\bibitem{QCQP1}
K.~M. Anstreicher.
\newblock Semidefinite programming versus the reformulation-linearization
  technique for nonconvex quadratically constrained quadratic programming.
\newblock {\em J. of Global Optimization}, 43(2):471--484, 2009.

\bibitem{mosek}
M.~ApS.
\newblock Mosek optimization toolbox for matlab.
\newblock {\em User’s Guide and Reference Manual, Version}, 4, 2019.

\bibitem{bai2019proximal}
Y.~Bai, J.~Duchi, and S.~Mei.
\newblock Proximal algorithms for constrained composite optimization, with
  applications to solving low-rank sdps.
\newblock {\em arXiv preprint arXiv:1903.00184}, 2019.

\bibitem{syn}
A.~S. Bandeira, N.~Boumal, and V.~Voroninski.
\newblock On the low-rank approach for semidefinite programs arising in
  synchronization and community detection.
\newblock In {\em Conference on Learning Theory}, pages 361--382. PMLR, 2016.

\bibitem{BAI}
A.~I. Barvinok.
\newblock Problems of distance geometry and convex properties of quadratic
  maps.
\newblock {\em Discrete \& Computational Geometry}, 13(2):189--202, 1995.

\bibitem{bellavia2021relaxed}
S.~Bellavia, J.~Gondzio, and M.~Porcelli.
\newblock A relaxed interior point method for low-rank semidefinite programming
  problems with applications to matrix completion.
\newblock {\em Journal of Scientific Computing}, 89(2):46, 2021.

\bibitem{Boumal1}
S.~Bhojanapalli, N.~Boumal, P.~Jain, and P.~Netrapalli.
\newblock Smoothed analysis for low-rank solutions to semidefinite programs in
  quadratic penalty form.
\newblock In {\em Conference On Learning Theory}, pages 3243--3270. PMLR, 2018.

\bibitem{SNL}
P.~Biswas, T.-C. Liang, K.-C. Toh, Y.~Ye, and T.-C. Wang.
\newblock Semidefinite programming approaches for sensor network localization
  with noisy distance measurements.
\newblock {\em IEEE Transactions on Automation Science and Engineering},
  3(4):360--371, 2006.

\bibitem{bonnans2013perturbation}
J.~F. Bonnans and A.~Shapiro.
\newblock {\em Perturbation analysis of optimization problems}.
\newblock Springer Science \& Business Media, 2013.

\bibitem{Staircase1}
N.~Boumal.
\newblock A {Riemannian} low-rank method for optimization over semidefinite
  matrices with block-diagonal constraints.
\newblock {\em arXiv preprint arXiv:1506.00575}, 2015.

\bibitem{Intromani}
N.~Boumal.
\newblock An introduction to optimization on smooth manifolds.
\newblock {\em Available online, May}, 3, 2020.

\bibitem{BouRTR}
N.~Boumal, P.-A. Absil, and C.~Cartis.
\newblock Global rates of convergence for nonconvex optimization on manifolds.
\newblock {\em IMA J. of Numerical Analysis}, 39(1):1--33, 2019.

\bibitem{manopt}
N.~Boumal, B.~Mishra, P.-A. Absil, and R.~Sepulchre.
\newblock Manopt, a {Matlab} toolbox for optimization on manifolds.
\newblock {\em J. of Machine Learning Research}, 15(1):1455--1459, 2014.

\bibitem{Boumal2}
N.~Boumal, V.~Voroninski, and A.~Bandeira.
\newblock The non-convex {Burer-Monteiro} approach works on smooth semidefinite
  programs.
\newblock {\em Advances in Neural Information Processing Systems}, 29, 2016.

\bibitem{Boumal3}
N.~Boumal, V.~Voroninski, and A.~S. Bandeira.
\newblock Deterministic guarantees for {Burer-Monteiro} factorizations of
  smooth semidefinite programs.
\newblock {\em Communications on Pure and Applied Mathematics}, 73(3):581--608,
  2020.

\bibitem{BM3}
S.~Burer and C.~Choi.
\newblock Computational enhancements in low-rank semidefinite programming.
\newblock {\em Optimisation Methods and Software}, 21(3):493--512, 2006.

\bibitem{BM1}
S.~Burer and R.~D. Monteiro.
\newblock A nonlinear programming algorithm for solving semidefinite programs
  via low-rank factorization.
\newblock {\em Mathematical Programming}, 95(2):329--357, 2003.

\bibitem{BM2}
S.~Burer and R.~D. Monteiro.
\newblock Local minima and convergence in low-rank semidefinite programming.
\newblock {\em Mathematical Programming}, 103(3):427--444, 2005.

\bibitem{chan2008constraint}
Z.~X. Chan and D.~Sun.
\newblock Constraint nondegeneracy, strong regularity, and nonsingularity in
  semidefinite programming.
\newblock {\em SIAM J. Optimization}, 19(1):370--396, 2008.

\bibitem{CIFrank}
D.~Cifuentes.
\newblock On the {Burer--Monteiro} method for general semidefinite programs.
\newblock {\em Optimization Letters}, 15(6):2299--2309, 2021.

\bibitem{erdogdu2022convergence}
M.~A. Erdogdu, A.~Ozdaglar, P.~A. Parrilo, and N.~D. Vanli.
\newblock Convergence rate of block-coordinate maximization {Burer--Monteiro}
  method for solving large {SDPs}.
\newblock {\em Mathematical Programming}, 195(1-2):243--281, 2022.

\bibitem{FSL}
S.~Friedland and R.~Loewy.
\newblock Subspaces of symmetric matrices containing matrices with a multiple
  first eigenvalue.
\newblock {\em Pacific J. of Mathematics}, 62(2):389--399, 1976.

\bibitem{rankadap}
B.~Gao and P.-A. Absil.
\newblock A {Riemannian} rank-adaptive method for low-rank matrix completion.
\newblock {\em Computational Optimization and Applications}, 81(1):67--90,
  2022.

\bibitem{RBB1}
B.~Gao, N.~T. Son, P.-A. Absil, and T.~Stykel.
\newblock Riemannian optimization on the symplectic {Stiefel} manifold.
\newblock {\em SIAM J. Optimization}, 31(2):1546--1575, 2021.

\bibitem{ge2017no}
R.~Ge, C.~Jin, and Y.~Zheng.
\newblock No spurious local minima in nonconvex low rank problems: A unified
  geometric analysis.
\newblock In {\em International Conference on Machine Learning}, pages
  1233--1242. PMLR, 2017.

\bibitem{DT}
M.~W. Hirsch.
\newblock {\em Differential topology}, volume~33.
\newblock Springer Science \& Business Media, 2012.

\bibitem{non-mono}
B.~Iannazzo and M.~Porcelli.
\newblock The {Riemannian Barzilai--Borwein} method with nonmonotone line
  search and the matrix geometric mean computation.
\newblock {\em IMA J. of Numerical Analysis}, 38(1):495--517, 2018.

\bibitem{Staircase2}
M.~Journ{\'e}e, F.~Bach, P.-A. Absil, and R.~Sepulchre.
\newblock Low-rank optimization on the cone of positive semidefinite matrices.
\newblock {\em SIAM J. Optimization}, 20(5):2327--2351, 2010.

\bibitem{karimi2016linear}
H.~Karimi, J.~Nutini, and M.~Schmidt.
\newblock Linear convergence of gradient and proximal-gradient methods under
  the {Polyak-{L}ojasiewicz} condition.
\newblock In {\em Machine Learning and Knowledge Discovery in Databases:
  European Conference, ECML PKDD 2016, Riva del Garda, Italy, September 19-23,
  2016, Proceedings, Part I 16}, pages 795--811. Springer, 2016.

\bibitem{Geoopt}
M.~Kochurov, R.~Karimov, and S.~Kozlukov.
\newblock Geoopt: {Riemannian} optimization in pytorch.
\newblock {\em arXiv preprint arXiv:2005.02819}, 2020.

\bibitem{lee2022escaping}
C.-p. Lee, L.~Liang, T.~Tang, and K.-C. Toh.
\newblock Accelerating nuclear-norm regularized low-rank matrix optimization
  through {Burer-Monteiro} decomposition.
\newblock {\em arXiv preprint arXiv:2204.14067}, 2022.

\bibitem{QSDPNAL}
X.~Li, D.~Sun, and K.-C. Toh.
\newblock {QSDPNAL:} a two-phase augmented {Lagrangian} method for convex
  quadratic semidefinite programming.
\newblock {\em Mathematical Programming Computation}, 10(4):703--743, 2018.

\bibitem{lojasiewicz1963topological}
S.~Lojasiewicz.
\newblock A topological property of real analytic subsets.
\newblock {\em Coll. du CNRS, Les {\'e}quations aux d{\'e}riv{\'e}es
  partielles}, 117(87-89):2, 1963.

\bibitem{lourencco2018optimality}
B.~F. Louren{\c{c}}o, E.~H. Fukuda, and M.~Fukushima.
\newblock Optimality conditions for nonlinear semidefinite programming via
  squared slack variables.
\newblock {\em Mathematical Programming}, 168:177--200, 2018.

\bibitem{Theta}
L.~Lov{\'a}sz.
\newblock On the {Shannon} capacity of a graph.
\newblock {\em IEEE Transactions on Information theory}, 25(1):1--7, 1979.

\bibitem{RegSDP}
J.~Malick, J.~Povh, F.~Rendl, and A.~Wiegele.
\newblock Regularization methods for semidefinite programming.
\newblock {\em SIAM J. Optimization}, 20(1):336--356, 2009.

\bibitem{PG}
G.~Pataki.
\newblock On the rank of extreme matrices in semidefinite programs and the
  multiplicity of optimal eigenvalues.
\newblock {\em Mathematics of Operations Research}, 23(2):339--358, 1998.

\bibitem{Boumal4}
T.~Pumir, S.~Jelassi, and N.~Boumal.
\newblock Smoothed analysis of the low-rank approach for smooth semidefinite
  programs.
\newblock {\em Advances in Neural Information Processing Systems}, 31, 2018.

\bibitem{Sardap}
S.~Scholtes and M.~St{\"o}hr.
\newblock How stringent is the linear independence assumption for mathematical
  programs with complementarity constraints?
\newblock {\em Mathematics of Operations Research}, 26(4):851--863, 2001.

\bibitem{Coding}
N.~J.~A. Sloane.
\newblock Challenge problems: Independent sets in graphs.
\newblock {\em http://www. research. att. com/\~{} njas/doc/graphs. html},
  2005.

\bibitem{SDPNALp2}
D.~Sun, K.-C. Toh, Y.~Yuan, and X.-Y. Zhao.
\newblock {SDPNAL+}: A {Matlab} software for semidefinite programming with
  bound constraints (version 1.0).
\newblock {\em Optimization Methods and Software}, 35(1):87--115, 2020.

\bibitem{GEP}
T.~Tang and K.-C. Toh.
\newblock Solving graph equipartition {SDPs} on an algebraic variety.
\newblock {\em arXiv preprint arXiv:2112.04256}, 2021.

\bibitem{SQK}
T.~Tang and K.-C. Toh.
\newblock A feasible method for solving an {SDP} relaxation of the quadratic
  knapsack problem.
\newblock {\em Mathematics of Operations Research, to appear}, 2023.

\bibitem{TTT}
K.-C. Toh, M.~J. Todd, and R.~H. T{\"u}t{\"u}nc{\"u}.
\newblock {SDPT3—a Matlab} software package for semidefinite programming,
  version 1.3.
\newblock {\em Optimization Methods and Software}, 11(1-4):545--581, 1999.

\bibitem{T3Q}
R.~H. T{\"u}t{\"u}nc{\"u}, K.-C. Toh, and M.~J. Todd.
\newblock Solving semidefinite-quadratic-linear programs using {SDPT3}.
\newblock {\em Mathematical Programming}, 95(2):189--217, 2003.

\bibitem{rankada1}
A.~Uschmajew and B.~Vandereycken.
\newblock Greedy rank updates combined with {Riemannian} descent methods for
  low-rank optimization.
\newblock In {\em International Conference on Sampling Theory and Applications
  (SampTA)}, pages 420--424. IEEE, 2015.

\bibitem{BQP2}
D.~Vandenbussche and G.~L. Nemhauser.
\newblock A branch-and-cut algorithm for nonconvex quadratic programs with box
  constraints.
\newblock {\em Mathematical Programming}, 102(3):559--575, 2005.

\bibitem{wang2023solving}
J.~Wang and L.~Hu.
\newblock Solving low-rank semidefinite programs via manifold optimization.
\newblock {\em arXiv preprint arXiv:2303.01722v1}, 2023.

\bibitem{wang2021decomposition}
Y.~Wang, K.~Deng, H.~Liu, and Z.~Wen.
\newblock A decomposition augmented lagrangian method for low-rank semidefinite
  programming.
\newblock {\em SIAM Journal on Optimization}, 33(3):1361--1390, 2023.

\bibitem{WYmani}
Z.~Wen and W.~Yin.
\newblock A feasible method for optimization with orthogonality constraints.
\newblock {\em Mathematical Programming}, 142(1):397--434, 2013.

\bibitem{CDopt1}
N.~Xiao, X.~Hu, X.~Liu, and K.-C. Toh.
\newblock Cdopt: A python package for a class of {Riemannian} optimization.
\newblock {\em arXiv preprint arXiv:2212.02698}, 2022.

\bibitem{CDopt}
N.~Xiao, X.~Liu, and K.-C. Toh.
\newblock Constraint dissolving approaches for {Riemannian} optimization.
\newblock {\em arXiv preprint arXiv:2203.10319}, 2022.

\bibitem{SDPNALp1}
L.~Yang, D.~Sun, and K.-C. Toh.
\newblock {SDPNAL+}: a majorized semismooth {Newton-CG augmented Lagrangian}
  method for semidefinite programming with nonnegative constraints.
\newblock {\em Mathematical Programming Computation}, 7(3):331--366, 2015.

\bibitem{CGAL}
A.~Yurtsever, J.~A. Tropp, O.~Fercoq, M.~Udell, and V.~Cevher.
\newblock Scalable semidefinite programming.
\newblock {\em SIAM J. on Mathematics of Data Science}, 3(1):171--200, 2021.

\bibitem{zhang2021preconditioned}
J.~Zhang, S.~Fattahi, and R.~Y. Zhang.
\newblock Preconditioned gradient descent for over-parameterized nonconvex
  matrix factorization.
\newblock {\em Advances in Neural Information Processing Systems},
  34:5985--5996, 2021.

\bibitem{newretrac}
R.~Zhang.
\newblock Newton retraction as approximate geodesics on submanifolds.
\newblock {\em arXiv preprint arXiv:2006.14751}, 2020.

\bibitem{zhang2022improved}
R.~Y. Zhang.
\newblock Improved global guarantees for the nonconvex burer--monteiro
  factorization via rank overparameterization.
\newblock {\em arXiv preprint arXiv:2207.01789}, 2022.

\bibitem{SDPNAL}
X.-Y. Zhao, D.~Sun, and K.-C. Toh.
\newblock A {Newton-CG augmented Lagrangian} method for semidefinite
  programming.
\newblock {\em SIAM J. Optimization}, 20(4):1737--1765, 2010.

\bibitem{rankada2}
G.~Zhou, W.~Huang, K.~A. Gallivan, P.~Van~Dooren, and P.-A. Absil.
\newblock A {Riemannian} rank-adaptive method for low-rank optimization.
\newblock {\em Neurocomputing}, 192:72--80, 2016.

\end{thebibliography}

\end{document}